\documentclass[11pt]{article}
\usepackage[utf8]{inputenc}
\usepackage[T1]{fontenc}
\usepackage[english]{babel}
\usepackage[letterpaper,margin=1in]{geometry}
\usepackage{amsmath,amssymb,amsthm,mathtools}
\usepackage{enumitem}
\usepackage{booktabs}
\usepackage{graphicx}
\usepackage{xcolor}
\usepackage{authblk}
\usepackage[mathlines]{lineno}
\usepackage[colorlinks=true,linkcolor=blue!60!black,citecolor=green!45!black,urlcolor=blue]{hyperref}
\usepackage{cleveref}
\usepackage[numbers, sort&compress]{natbib}
\usepackage[font=small,labelfont=bf]{caption}

\newtheorem{teorema}{Theorem}[section]
\newtheorem{lema}[teorema]{Lemma}
\newtheorem{proposicion}[teorema]{Proposition}
\newtheorem{corolario}[teorema]{Corollary}
\newtheorem{conjetura}[teorema]{Conjecture}
\newtheorem{problema}[teorema]{Problem}
\theoremstyle{definition}

\newtheorem{ejemplo}[teorema]{Example}
\newtheorem{construccion}[teorema]{Construction}
\theoremstyle{remark}

\theoremstyle{plain}
\newtheorem*{teoremarepinner}{Theorem \theteoremarepnum}
\newcommand{\theteoremarepnum}{}
\newenvironment{teoremarep}[1]{\renewcommand{\theteoremarepnum}{\ref{r:#1}}\begin{teoremarepinner}}{\end{teoremarepinner}}

\newcommand{\cs}[1]{\mathcal{#1}}
\DeclareMathOperator{\cro}{cr}
\DeclareMathOperator{\gen}{gen}
\DeclareMathOperator{\Sym}{Sym}
\DeclareMathOperator{\Aut}{Aut}
\DeclareMathOperator{\ex}{ex}

\usepackage{environ}

\title{Token graphs of stars}
\author[1]{Emmanuel Acosta Troncoso\thanks{emmanuel\_acosta22@yahoo.com}}
\author[2]{César Hernández-Cruz\thanks{Supported by UNAM-PAPIIT grant IN106425 and SEP-CONACYT Ciencia Básica A1-S-8397 grant. \texttt{chc@ciencias.unam.mx}}}
\author[2,3]{Ana Trujillo-Negrete\thanks{Supported by the Universidad Nacional Autónoma de México Postdoctoral Program (POSDOC). \texttt{atrujillo@im.unam.mx}}}
\affil[1]{\textit{Unidad Académica de Matemáticas, Universidad Autónoma de Zacatecas, México}}
\affil[2]{\textit{Facultad de Ciencias, Universidad Nacional Autónoma de México, México}}
\affil[3]{\textit{Instituto de Matemáticas, Universidad Nacional Autónoma de México, México}}
\date{}

\begin{document}

\maketitle

\begin{abstract}
The $k$-token graph $F_k(G)$ of a graph $G$ has as vertices the $k$-subsets of $V(G)$, two of them being adjacent when their symmetric difference is an edge of $G$. Let $S_n$ be the star with $n$ leaves. Since every edge of $S_n$ is incident with the center, $F_k(S_n)$ is the inclusion graph of the $(k-1)$- and the $k$-subsets of the leaves, that is, the bipartite graph in which a $(k-1)$-subset is adjacent to the $k$-subsets that contain it. We use this identification to determine its diameter, radius, girth, Wiener index and number of spanning trees, and to decompose $F_k(S_{n+1})$ into $F_k(S_n)$ and $F_{k-1}(S_n)$ joined by a matching; the parameters of $F_k(S_n)$ that are already known, such as its connectivity, independence number, automorphism group and spectrum, are collected in an appendix. The two parameters that the identification does not settle are the domination and the packing numbers. For the domination number we prove, for fixed $k\ge 3$, that
\[
\Bigl(1-\tfrac{k-2}{k-1}\,\pi_k-o(1)\Bigr)\binom{n}{k-1}\;\le\;\gamma\bigl(F_k(S_n)\bigr)\;\le\;\Bigl(1-\tfrac12\bigl(\tfrac{k-2}{k-1}\bigr)^{k-2}+o(1)\Bigr)\binom{n}{k-1},
\]
where $\pi_k$ is the Turán density of the complete $(k-1)$-uniform hypergraph on $k$ vertices; for $k=3$ both constants are equal to $\frac34$, so that $\gamma(F_3(S_n))=\bigl(\frac34+o(1)\bigr)\binom n2$, and the upper bound takes the exact form $\gamma(F_3(S_n))\le\binom n2-\lfloor (n+1)^2/8\rfloor$, which is attained for all $3\le n\le 10$. For $k=4$ we improve the upper bound to $\bigl(\frac{17}{27}+o(1)\bigr)\binom{n}{3}$, which matches the lower bound if Turán's conjecture $\pi_4=\frac59$ holds. We show that the packing number of $F_k(S_n)$ equals the maximum size $A(n+1,4,k)$ of a binary constant-weight code, and we characterize when the packing and the domination numbers coincide, namely, for $2\le k\le n-1$, exactly when $n=2k-1$ and a Steiner system $S(k-2,k-1,2k-1)$ exists; in that case both parameters equal the Catalan number $C_k$.
\end{abstract}

{\small\noindent\textbf{Keywords:} token graphs, star, domination number, packing number, Turán density, Steiner systems.

\noindent\textbf{MSC 2020:} 05C69, 05C65, 05C35, 05B40.}

\section{Introduction}\label{sec:1}

\begin{figure}[t]
    \centering
    \includegraphics[width=0.75\linewidth]{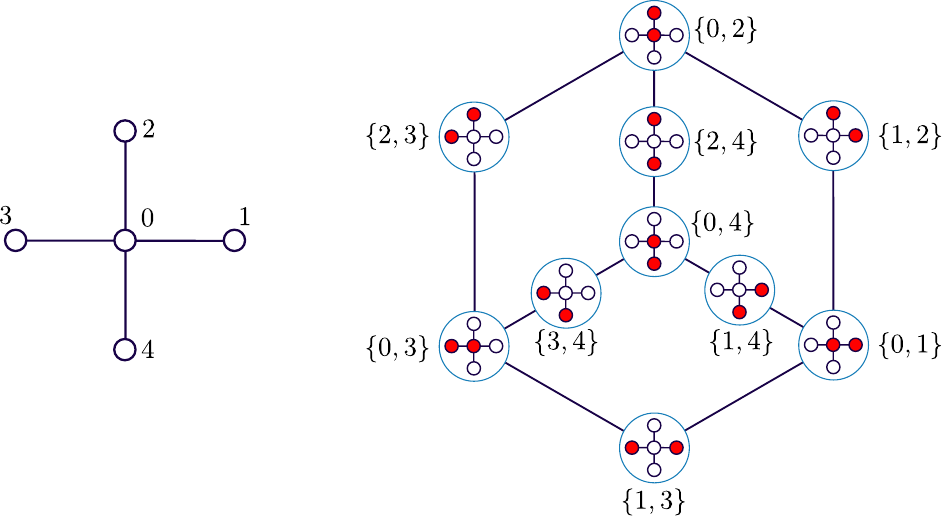}
    \caption{The star $S_4$ on the left and its $2$-token graph $F_2(S_4)$ on the right; each vertex is drawn as the placement of tokens it represents.}
    \label{fig:example}
\end{figure}

Let $G$ be a simple graph and let $k$ be an integer with $1\le k\le |V(G)|-1$. The \emph{$k$-token graph} $F_k(G)$ has as vertices the $k$-subsets of $V(G)$, two of them being adjacent when their symmetric difference is an edge of $G$; that is, $A$ and $B$ are adjacent when there are $u\in A$ and $v\in B$ with $uv\in E(G)$ and $A\setminus\{u\}=B\setminus\{v\}$. Equivalently, think of $k$ indistinguishable tokens occupying distinct vertices of $G$. A vertex of $F_k(G)$ is then a placement of the tokens, and an edge is the sliding of one token along an edge of $G$ to an empty vertex. Thus $F_k(G)$ is the configuration space of $k$ tokens on $G$, and paths in $F_k(G)$ are exactly the valid sequences of moves. See an example in Figure~\ref{fig:example}.

These graphs were introduced under this name by Fabila-Monroy, Flores-Peñaloza, Huemer, Hurtado, Urrutia, and Wood~\cite{b25}, who established their basic properties and posed several of the questions that have guided their subsequent study. They had appeared before under other names: as \emph{$k$-subgraph graphs}~\cite{b36}, as \emph{double vertex graphs} ($k=2$) and \emph{$k$-tuple vertex graphs}~\cite{b3,b56}, and as \emph{symmetric $k$-th powers}~\cite{b51,b5}, the latter motivated by quantum mechanics. For $G=K_n$ one obtains the Johnson graph $J(n,k)$, a fact that will reappear twice below. 
Among the parameters of $F_k(G)$ that have been studied are connectivity~\cite{b25,b42,b26}, planarity and regularity~\cite{b14}, independence and matching numbers~\cite{b19}, Laplacian spectra and algebraic connectivity~\cite{b18,b17}, automorphisms~\cite{b34,b24} and treewidth~\cite{b27}. Most recently, the domination number has been studied for the token graphs of paths~\cite{bALR,bAR}, of complete bipartite graphs~\cite{bCCPSS,bACR} and of the star~\cite{b23}; we compare our bounds with those of~\cite{b23} in Subsection~\ref{sec:4.5}.

In this work we take $G$ to be a star. We write $S_n=K_{1,n}$ for the star with $n$ leaves and fix once and for all the labeling
\[
V(S_n)=\{0,1,\dots,n\},\qquad E(S_n)=\{\,0\, i : 1\le i\le n\,\},
\]
so that $0$ is the center, $1,\dots,n$ are the leaves, and $S_n$ has $n+1$ vertices. Every edge of $S_n$ is incident with the center, so a token can move only between the center and a leaf. Consequently a configuration is either \emph{centered}, if it contains $0$, or \emph{center-free}, if it does not, and no elementary move joins two configurations of the same kind. This dichotomy is the only structural fact about $S_n$ that we shall use, and it organizes everything that follows. Forgetting the center, a centered configuration is a $(k-1)$-subset $X$ of the leaves and a center-free one is a $k$-subset $A$, and a move joins $X$ to $A$ exactly when $X\subset A$. In other words, $F_k(S_n)$ is the \emph{inclusion graph} $B(n;k-1,k)$, the bipartite graph whose vertices are the $(k-1)$- and the $k$-subsets of $[n]$, a $(k-1)$-subset being adjacent to the $k$-subsets that contain it (Theorem~\ref{r:3.1}). This graph is the \emph{doubled Johnson graph} of Dalfó et al.~\cite{b18}, and for $n=2k-1$ it is the \emph{middle levels graph}, whose Hamiltonicity was established by Mütze~\cite{b43}.

The aforementioned identification is the engine of the paper, and it is remarkably efficient. The distance between two vertices of $B(n;k-1,k)$ is determined by the number of elements in which the two subsets differ (Lemma~\ref{r:3.6}), which gives the diameter and the radius. The symmetric group $\Sym(n)$ acts edge-transitively on $B(n;k-1,k)$, which gives the connectivity through a theorem of Watkins.   Since $B(n;k-1,k)$ is bipartite and Hall's condition holds, we obtain its independence, matching and vertex-cover numbers. As the Gram matrix of the inclusion matrix of $B(n;k-1,k)$ is a polynomial in the adjacency matrix of $J(n,k-1)$, the spectrum and the Laplacian spectrum are also obtained. Section~\ref{sec:3} and Appendix~\ref{sec:A} collect these facts, together with the exact values of the girth, chromatic number, chromatic index, Wiener index, and number of spanning trees.  Also in Section~\ref{sec:3}, we present an inductive decomposition (Theorem~\ref{r:3.3}),
\[
F_k(S_{n+1})\;\cong\;F_k(S_n)\;\sqcup\;F_{k-1}(S_n)\;+\;M,
\]
which realizes Pascal's rule $\binom{n+1}{k}=\binom{n}{k}+\binom{n}{k-1}$ as a decomposition of graphs, where $M$ is a matching between the centered configurations of $F_k(S_n)$ and the center-free configurations of $F_{k-1}(S_n)$.   Characterizations of the automorphism group for every graph in the family as well as planar members of the family are also presented, in Appendix~\ref{sec:A}.

Two natural parameters are not settled by the identification alone, namely the domination number $\gamma$ and the packing number $\rho$ (the largest number of vertices with pairwise disjoint closed neighborhoods). Both become problems about set systems, of two different kinds, and they occupy the rest of the paper. Packing is the easier to translate. A family of $r$-subsets of a $v$-set is a \emph{Steiner system} $S(t,r,v)$ if every $t$-subset of the $v$-set is contained in exactly one of its members, and a \emph{partial} $S(t,r,v)$ if every $t$-subset is contained in at most one of them; we write $\mathcal{D}(t,r,v)$ for the maximum size of a partial $S(t,r,v)$. The case of interest here is $t=k-1$ and $r=k$. A family of $k$-subsets is a partial $S(k-1,k,v)$ if and only if any two of its members share at most $k-2$ points, and $\mathcal{D}(k-1,k,v)$ is also the maximum size $A(v,4,k)$ of a binary code of length $v$, constant weight $k$ and minimum distance $4$. Since two configurations of $F_k(S_n)$ are at distance at most $2$ exactly when they share $k-1$ points, we obtain the following connection with design and coding theory.

\begin{teorema}\label{r:1.1}
Let $1\le k\le n$. A set $P$ of vertices of $F_k(S_n)$ is a packing if and only if $|A\cap B|\le k-2$ for all $A\ne B$ in $P$; that is, if and only if $P$ is a partial $S(k-1,k,n+1)$ on $V(S_n)$. Consequently,
\[
\rho(F_k(S_n))=\mathcal{D}(k-1,k,n+1)=A(n+1,4,k).
\]
\end{teorema}

Thus the packing number of $F_k(S_n)$ is a classical quantity of design and coding theory, and the results on $\mathcal{D}(k-1,k,v)$ transfer to it, in particular the Johnson bound~\cite{b37}, Rödl's asymptotic formula~\cite{b50}, Keevash's existence theorem for designs~\cite{b38}, and the exact values for $k=3$~\cite{b39,b52,b53} and $k=4$~\cite{b32,b9,b7}. Rödl's and Keevash's theorems are discussed in Section~\ref{sec:5}, together with two explicit constructions by means of abelian groups; the Johnson bound and the exact values are collected in Appendix~\ref{sec:B}.

Domination does not translate as directly. A dominating set $D$ of $B(n;k-1,k)$ is best described by the $(k-1)$-sets that it does \emph{not} contain, which form a $(k-1)$-uniform hypergraph $\mathcal{H}$, because $D$ must then contain every $k$-set all of whose $(k-1)$-subsets are edges of $\mathcal{H}$ (Proposition~\ref{r:4.1}). Minimizing $|D|$ therefore means finding hypergraphs with many edges and few such \emph{complete} $k$-sets, which is a problem of extremal hypergraph theory. A construction built from books gives the upper bound in Theorem~\ref{r:1.2} for every $k$, a hypergraph of Turán~\cite{b54} gives a better one for $k=4$, and the hypergraph removal lemma~\cite{bRSz,bNRS,bRS06,bGowers} gives the lower bound, in terms of the Turán density $\pi_k$ of the complete $(k-1)$-uniform hypergraph $K_k^{(k-1)}$ on $k$ vertices. For $k=2$ the value is known exactly, $\gamma(F_2(S_n))=n-1$~\cite{b23}.

\begin{teorema}\label{r:1.2}
Let $k\ge 3$ be fixed. Then, as $n\to\infty$,
\[
\Bigl(1-\tfrac{k-2}{k-1}\,\pi_k-o(1)\Bigr)\binom{n}{k-1}\;\le\;\gamma\bigl(F_k(S_n)\bigr)\;\le\;\Bigl(1-\tfrac12\bigl(\tfrac{k-2}{k-1}\bigr)^{k-2}+o(1)\Bigr)\binom{n}{k-1}.
\]
In particular, $\gamma(F_k(S_n))\ge\bigl(\frac{2k-3}{(k-1)^2}-o(1)\bigr)\binom{n}{k-1}$, and for $k=3$,
\[
\bigl(\tfrac34-o(1)\bigr)\binom{n}{2}\;\le\;\gamma\bigl(F_3(S_n)\bigr)\;\le\;\binom{n}{2}-\Bigl\lfloor\frac{(n+1)^2}{8}\Bigr\rfloor\;=\;\bigl(\tfrac34+o(1)\bigr)\binom{n}{2};
\]
that is, $\gamma(F_3(S_n))=\bigl(\frac34+o(1)\bigr)\binom n2$. Moreover, $\gamma(F_4(S_n))\le\bigl(\frac{17}{27}+o(1)\bigr)\binom{n}{3}$. 
\end{teorema}

The lower bound more than doubles the bound $\frac1k\binom{n}{k-1}$ of Fabila-Monroy and Gómez-Galicia~\cite{b23}, and the upper bound is better than theirs when $k=4$, when $k\in\{6,26,36,56\}$, and whenever the smallest prime factor of $k-1$ is at least $7$ (Subsection~\ref{sec:4.5}). For $k=4$, the bound $\pi_4\le0.5578$ of~\cite{b35} gives the lower constant $0.628$, against the upper constant $\frac{17}{27}\approx 0.6296$, and the two coincide if Turán's conjecture $\pi_4=\frac59$ holds (Theorem~\ref{r:4.12}). The exact upper bound for $k=3$ is attained for all $3\le n\le 10$, is asymptotically correct by Theorem~\ref{r:1.2}, and we conjecture it to be the right value (Conjecture~\ref{r:7.2}).

Domination and packing meet in Section~\ref{sec:6}. The Steiner bound for packing and the volume bound for domination are the same number, $\frac1k\binom{n+1}{k-1}=|V(F_k(S_n))|/(\Delta(F_k(S_n))+1)$, and $\rho=\gamma$ forces this number to be attained by a perfect dominating set consisting of vertices of maximum degree. For $2\le k\le n-1$ the existence of such a set is possible only in the regular case $n=2k-1$, where the perfect dominating sets of $F_k(S_{2k-1})$ correspond exactly to the Steiner systems $S(k-2,k-1,2k-1)$ (Theorem~\ref{r:6.3}).  Since $F_k(S_{2k-1})$ is the bipartite double of the odd graph $O_k$, and a perfect dominating set of the former corresponds to a perfect $1$-code of the latter, this is also a theorem of Hammond and Smith~\cite{b31} on perfect codes in odd graphs.  For $k\in\{1,n\}$ the graph $F_k(S_n)$ is the star itself and $\rho=\gamma=1$, so these trivial cases are left out. Therefore, the following result is a complete characterization.

\begin{teorema}\label{r:1.3}
Let $2\le k\le n-1$. Then $\rho(F_k(S_n))=\gamma(F_k(S_n))$ if and only if $n=2k-1$ and there exists a Steiner system $S(k-2,k-1,2k-1)$, in which case $\rho(F_k(S_n))=\gamma(F_k(S_n))=C_k=\frac{1}{k+1}\binom{2k}{k}$, the $k$-th Catalan number.
\end{teorema}

Since a perfect dominating set is at the same time a packing and a dominating set, Theorem~\ref{r:1.3} also says that, for $2\le k\le n-1$, the graph $F_k(S_n)$ has a perfect dominating set only when $n=2k-1$ (Corollary~\ref{r:6.8}). Equality occurs for $k\in\{2,4,6\}$, through the trivial system $S(0,1,3)$, the Fano plane $S(2,3,7)$ and the Witt design $S(4,5,11)$, and it requires $k+1$ to be prime; the first undecided case is $k=10$.

The rest of the paper is organized as follows. Section~\ref{sec:2} contains the preliminaries and fixes notation. Section~\ref{sec:3} proves the identification of $F_k(S_n)$ and $B(n;k-1,k)$, presents the inductive decomposition, and computes the distances, the diameter, the radius, the Wiener index and the number of spanning trees, and settles Hamiltonicity. Section~\ref{sec:4} treats domination, namely the hypergraph reformulation, the book construction, the removal-lemma lower bound, the cases $k\le 4$, and the comparison with~\cite{b23}. Section~\ref{sec:5} treats packing: it proves Theorem~\ref{r:1.1}, derives the Steiner bound and the asymptotic value of $\rho(F_k(S_n))$, illustrates the correspondence with designs on small cases, and gives two constructions of packings by means of abelian groups. Section~\ref{sec:6} proves Theorem~\ref{r:1.3}, and Section~\ref{sec:7} collects some open problems. The appendices contain what is known or routine. In Appendix~\ref{sec:A} we prove, through $B(n;k-1,k)$, several parameters of $F_k(S_n)$ that were already known; in Appendix~\ref{sec:B} we state the results from design and coding theory that we use; and in Appendix~\ref{sec:C} we give two proofs that were deferred from the text, together with examples of these constructions and their comparison with the exact value for $k=3$.

\section{Preliminaries}\label{sec:2}

We follow~\cite{b21} for standard graph-theoretic notation. All graphs are finite and simple, and $G$ always denotes a graph. We write $[n]=\{1,\dots,n\}$, $\binom{X}{r}$ for the family of $r$-subsets of a set $X$, and $\Sym(n)$ for the symmetric group on $[n]$.

\medskip\noindent\textit{Degrees and structure.} The minimum and maximum degree of $G$ are denoted by $\delta(G)$ and $\Delta(G)$, respectively. A bipartite graph is \emph{biregular} if the vertices on each side of the bipartition have the same degree. The \emph{total subdivision} of $G$ is the graph obtained by subdividing every edge of $G$ exactly once, so that each edge becomes a path of length two through a new vertex of degree two. We write $\Aut(G)$ for the automorphism group of $G$, that is, the group of permutations of $V(G)$ preserving adjacency.

\medskip\noindent\textit{Distances.} The distance between $u$ and $v$ is denoted $d(u,v)$, the eccentricity of $u$ is $\max_v d(u,v)$, and $\operatorname{diam}(G)$ and $\operatorname{rad}(G)$ are the maximum and the minimum eccentricity; the \emph{center} of $G$ is the set of vertices of minimum eccentricity. The \emph{girth}, denoted $\operatorname{girth}(G)$, is the length of a shortest cycle, with $\operatorname{girth}(G)=\infty$ if $G$ is acyclic. The \emph{Wiener index} is the sum of all distances,
\[
W(G)=\sum_{\{u,v\}\subseteq V(G)} d(u,v).
\]

\medskip\noindent\textit{Colorings, independence, connectivity.} The chromatic number and the chromatic index are denoted by $\chi(G)$ and $\chi'(G)$; $\alpha(G)$, $\nu(G)$ and $\tau(G)$ are the independence, matching and vertex-cover numbers; and $\kappa(G)$ and $\lambda(G)$ are the vertex- and edge-connectivity.

\medskip\noindent\textit{Domination and packing.} A \emph{dominating set} of $G$ is a set $D\subseteq V(G)$ such that every vertex lies in $D$ or has a neighbor in $D$; the \emph{domination number} $\gamma(G)$ is the size of a smallest dominating set. The set $D$ is a \emph{perfect dominating set} (also \emph{efficient dominating set}, or \emph{perfect code}) if the closed neighborhoods $\{N[v]:v\in D\}$ partition $V(G)$, that is, if every vertex is dominated exactly once. A \emph{packing} is a set of vertices at pairwise distance at least $3$, equivalently a set whose closed neighborhoods are pairwise disjoint; the \emph{packing number} $\rho(G)$ is the largest size of a packing. In every graph $\rho(G)\le\gamma(G)$, since a dominating set must contain at least one vertex from each of the disjoint closed neighborhoods of a packing.

\medskip\noindent\textit{Steiner systems.} Let $t,r,v$ be integers with $0\le t<r<v$. A \emph{Steiner system} $S(t,r,v)$ is a family of $r$-subsets (called \emph{blocks}) of a $v$-set in which every $t$-subset is contained in exactly one block. A \emph{partial} $S(t,r,v)$ is one in which every $t$-subset is contained in at most one block, so that $t$-subsets may remain uncovered; equivalently, a family of $r$-subsets of a $v$-set is a partial $S(t,r,v)$ if and only if any two of its blocks meet in at most $t-1$ points. We write $\mathcal{D}(t,r,v)$ for the largest number of blocks of a partial $S(t,r,v)$,\footnote{Throughout the paper the parameters of Steiner systems, partial Steiner systems and their packing numbers are written in the same order $(t,r,v)$. The design-theory literature usually writes the packing number as $\mathcal{D}(v,r,t)$; we align it with $S(t,r,v)$ so that the two symbols can be read side by side without reversing the convention inside a single formula. For constant-weight codes we keep the standard notation $A(v,d,r)$ of Section~\ref{sec:5}.} and we allow $t=0$, where the only requirement is that there be exactly one block. The arithmetic of these systems is developed in Section~\ref{sec:6.2}.

\medskip\noindent\textit{Two reductions on the range of $k$.} Two elementary isomorphisms, both presented in~\cite{b25}, allow us to restrict the range of $k$ throughout. The first is $F_1(G)\cong G$. For the second, let $G$ be a graph on $n$ vertices and write $\mathfrak{c}(X)=V(G)\setminus X$ for the complement of a set of vertices. Then $\mathfrak{c}$ is a bijection from the $k$-subsets of $V(G)$ onto the $(n-k)$-subsets, and it leaves symmetric differences unchanged, $\mathfrak{c}(X)\mathbin{\triangle}\mathfrak{c}(Y)=X\mathbin{\triangle}Y$, so it preserves and reflects adjacency; hence $F_k(G)\cong F_{n-k}(G)$. Since $|V(S_n)|=n+1$, for the star this reads
\begin{equation}\label{eq:1}
F_k(S_n)\;\cong\;F_{n+1-k}(S_n),
\end{equation}
and every parameter of $F_k(S_n)$ is therefore invariant under $k\mapsto n+1-k$. Together with $F_1(S_n)\cong S_n$, this means that we may assume
\[
2\;\le\;k\;\le\;\tfrac{n+1}{2},
\]
which we do from Subsection~\ref{sec:3.3} onwards (the identification of Subsection~\ref{sec:3.1} and the inductive construction are stated for all $k$, because the construction uses $F_{k-1}(S_n)$). The extreme case $2k=n+1$ of this range is the \emph{self-complementary} one, where $F_k(S_n)\cong F_k(S_{2k-1})$ is regular; the sharpest results of Sections~\ref{sec:5} and~\ref{sec:6} concern precisely that case.

\section{Structure and parameters for arbitrary $k$}\label{sec:3}

Every edge of $S_n$ is incident with the center, so a token can move only between the center and a leaf; this naturally gives a bipartition of the token graph into those configurations that contain the center and those that do not. Moreover, as we show below in Theorem~\ref{r:3.1}, this bipartition already determines the graph and identifies it with a purely set-theoretic object, the containment graph between the $(k-1)$- and the $k$-subsets of $[n]$.

\subsection{The identification and an inductive construction}\label{sec:3.1}

The vertices of $F_k(S_n)$ split as mentioned above. There are $\binom{n}{k-1}$ \emph{centered} configurations, namely the sets $\{0\}\cup X$ with $X\in\binom{[n]}{k-1}$, and $\binom{n}{k}$ \emph{center-free} ones, namely the sets $A\in\binom{[n]}{k}$. Since every edge of $S_n$ has the form $0v$, for some $v\in[n]$, no two configurations of the same type can be adjacent, so the partition by type is a bipartition.

The \emph{inclusion graph} $B(n;k-1,k)$ is the bipartite graph with vertex set $\binom{[n]}{k-1}\cup\binom{[n]}{k}$, in which $X\in\binom{[n]}{k-1}$ and $A\in\binom{[n]}{k}$ are adjacent if and only if $X\subset A$. Throughout the paper, $X,Y,Z$ denote $(k-1)$-subsets of $[n]$ and $A,B,C$ denote $k$-subsets, while lowercase letters denote vertices of $S_n$. This graph appears as the \emph{doubled Johnson graph} $J(n;k-1,k)$ in Dalfó et al.~\cite[Section~5]{b18}, as a generalization of the doubled odd graphs of~\cite[Section~9.1D]{b11}, and its identification with the token graph of the star is observed there;\footnote{In~\cite{b18} the star $S_m$ has $m$ vertices, so our $F_k(S_n)$ is their $F_k(S_{n+1})$.} since everything in this paper depends on it, we include the short proof.

\begin{teorema}\label{r:3.1}
Let $1\le k\le n$. The map $\varphi\colon V(F_k(S_n))\to V\bigl(B(n;k-1,k)\bigr)$ given by
\[
\varphi(\{0\}\cup X)=X \qquad\text{and}\qquad \varphi(A)=A
\]
on the centered and the center-free configurations, respectively, is an isomorphism.
\end{teorema}

\begin{proof}
The map $\varphi$ is a bijection. Recall that two configurations are adjacent in $F_k(S_n)$ if and only if their symmetric difference is $\{0,v\}$ for some $v\in[n]$.

The symmetric difference of two centered configurations $\{0\}\cup X$ and $\{0\}\cup X'$ is $X\triangle X'$, which omits $0$. Therefore, these configurations are not adjacent, and their images lie on the same side of $B(n;k-1,k)$ and are not adjacent either. The same argument applies to two center-free configurations.

For a mixed pair, $\{0\}\cup X$ and $A$ are adjacent if and only if their symmetric difference is $\{0,v\}$ for some $v\in[n]$, that is, if and only if $A=X\cup\{v\}$ with $v\notin X$, that is, if and only if $X\subset A$. This is precisely the adjacency rule of $B(n;k-1,k)$.
\end{proof}

\begin{figure}[t]
\centering
\includegraphics[width=0.85\textwidth]{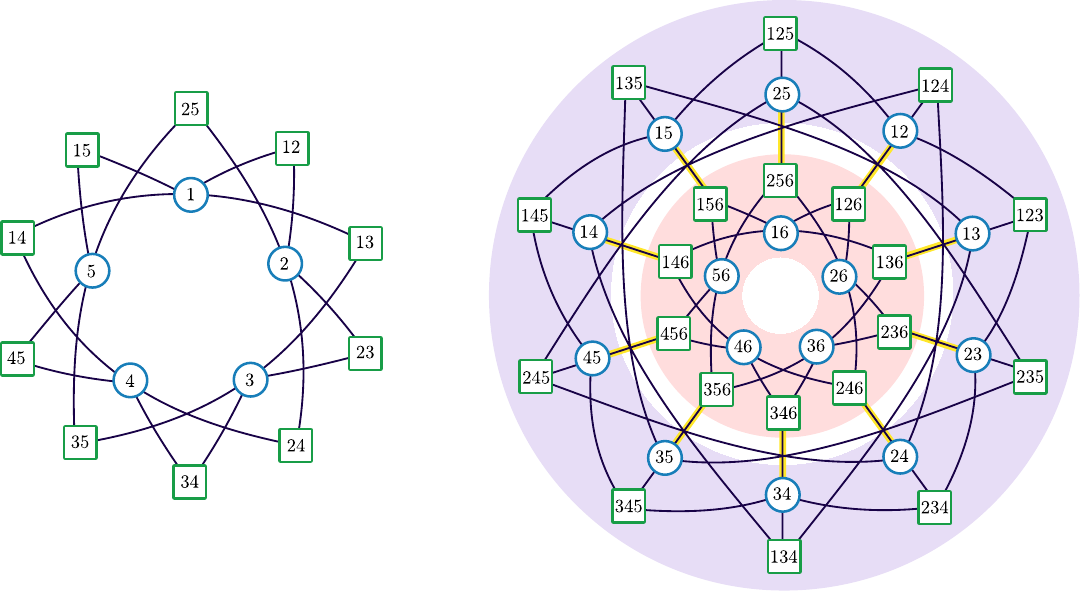}
\caption{The inclusion graphs $B(5;1,2)\cong F_2(S_5)$ (left) and $B(6;2,3)\cong F_3(S_6)$ (right). The $(k-1)$-subsets are drawn as blue circles and the $k$-subsets as green squares, and two of them are adjacent when one contains the other. On the right, the subsets that do not contain $6$ (purple region) induce $B(5;2,3)\cong F_3(S_5)$, those that contain $6$ (pink region) induce a copy of the graph on the left, and the yellow edges form the matching $M$ of Theorem~\ref{r:3.3}, which joins each $2$-subset $X$ of $[5]$ to $X\cup\{6\}$.}
\label{fig:1}
\end{figure}

The graphs $F_2(S_5)$ and $F_3(S_6)$ are depicted in Figure~\ref{fig:1} as their corresponding inclusion graphs.

\begin{corolario}\label{r:3.2}
For $1\le k\le n$ the graph $F_k(S_n)$ is bipartite and biregular, with
\[
\bigl|V(F_k(S_n))\bigr|=\binom{n}{k-1}+\binom{n}{k}=\binom{n+1}{k},\qquad \bigl|E(F_k(S_n))\bigr|=k\binom{n}{k}=(n-k+1)\binom{n}{k-1};
\]
the center-free configurations have degree $k$ and the centered ones have degree $n-k+1$; thus, $\delta(F_k(S_n))=\min\{k,n-k+1\}$ and $\Delta(F_k(S_n))=\max\{k,n-k+1\}$. 
In particular $F_k(S_n)$ is regular if and only if $2k=n+1$.
\end{corolario}

\begin{proof}
A $k$-subset $A$ contains exactly $k$ subsets of size $k-1$, and a $(k-1)$-subset $X$ extends to a $k$-subset in exactly $n-k+1$ ways; double counting the edges gives $|E(F_k(S_n))|$, and Pascal's rule gives $|V(F_k(S_n))|$. The bipartition is immediate from the definition of $B(n;k-1,k)$, and the two degrees coincide exactly when $k=n-k+1$, that is, when $2k=n+1$.
\end{proof}

We close this subsection with an inductive construction of $F_k(S_{n+1})$ in terms of $F_k(S_n)$ and $F_{k-1}(S_n)$.

\begin{teorema}\label{r:3.3}
For $2\le k\le n$,
\[
F_k(S_{n+1})\;\cong\;F_k(S_n)\;\sqcup\;F_{k-1}(S_n)\;+\;M,
\]
where the two token graphs appear as vertex-disjoint induced subgraphs and $M$ is the matching that joins, for each $(k-1)$-subset $X$ of $[n]$, the centered configuration $\{0\}\cup X$ of $F_k(S_n)$ to the center-free configuration $X$ of $F_{k-1}(S_n)$.
\end{teorema}

\begin{proof}
Let $z=n+1$ be the new leaf. The two copies are visible already in the language of tokens. A configuration of $F_k(S_{n+1})$ either leaves $z$ empty, and is then a configuration of $F_k(S_n)$, or keeps one token on $z$, and the remaining $k-1$ tokens form a configuration of $F_{k-1}(S_n)$. We verify that this splitting is the stated decomposition, and by Theorem~\ref{r:3.1} we argue throughout in the corresponding inclusion graphs, where a vertex of $B(n+1;k-1,k)$ is a $(k-1)$- or a $k$-subset of $[n+1]$, and the center no longer appears. Split the vertices according to whether they contain $z$,
\[
\cs{U}=\{\,S: z\notin S\,\},\qquad \cs{W}=\{\,S: z\in S\,\}.
\]

\emph{The copy $\cs{U}$.} A subset avoiding $z$ is a subset of $[n]$, so $\cs{U}$ consists of the $(k-1)$- and the $k$-subsets of $[n]$ with the induced containment relation, which is $B(n;k-1,k)\cong F_k(S_n)$.

\emph{The copy $\cs{W}$.} A subset containing $z$ is $S'\cup\{z\}$ with $S'\subseteq[n]$; removing $z$ sends the $(k-1)$-subsets to $(k-2)$-subsets and the $k$-subsets to $(k-1)$-subsets of $[n]$, and preserves containment, since $S'\cup\{z\}\subset T'\cup\{z\}$ if and only if $S'\subset T'$. Hence $\cs{W}$ induces $B(n;k-2,k-1)\cong F_{k-1}(S_n)$.

\emph{The cross edges.} An edge between $\cs{U}$ and $\cs{W}$ joins $X\in\cs{U}$ to $A\in\cs{W}$ with $X\subset A$; by cardinality $X$ is the $(k-1)$-set and $A$ the $k$-set. Writing $A=A'\cup\{z\}$ with $A'\subseteq[n]$, the containment $X\subset A$ forces $X\subseteq A'$, and $|X|=|A'|=k-1$ gives $X=A'$. So each $(k-1)$-subset $X$ of $[n]$ is joined to the unique vertex $X\cup\{z\}$, and there are no other cross edges. Under the two identifications $X$ is a centered vertex of the first copy and $X\cup\{z\}$ is a center-free vertex of the second, so these $\binom{n}{k-1}$ edges form the stated matching.
\end{proof}

Figure~\ref{fig:1} shows the case $k=3$, $n=5$ of the decomposition. In the graph on the right the fifteen vertices that contain the leaf $6$ induce a copy of the graph on the left, which is $B(5;1,2)\cong F_2(S_5)$, and the twenty remaining ones induce $B(5;2,3)\cong F_3(S_5)$. The matching $M$ is visible there as well, since it joins each $2$-subset $X$ of $[5]$, which is a vertex of the copy of $F_3(S_5)$, to the $3$-subset $X\cup\{6\}$, which under the identification of the copy on the left is the center-free configuration $X$ of $F_2(S_5)$. Iterating Theorem~\ref{r:3.3} expresses $F_k(S_n)$ through the triangle of token graphs $F_j(S_m)$ with $j\le k$ and $m<n$.

For $k=2$, Theorem~\ref{r:3.1} identifies $F_2(S_n)$ with the total subdivision of $K_n$. The vertices of $B(n;1,2)$ are the $1$- and the $2$-subsets of $[n]$, joined by containment; a center-free configuration $\{i,j\}$ is adjacent exactly to $\{0,i\}$ and $\{0,j\}$, so it has degree $2$ and its two neighbors are centered, while a centered configuration $\{0,i\}$ has degree $n-1$. Every center-free vertex therefore plays the role of the vertex that subdivides an edge.

\begin{corolario}\label{r:3.4}
For $n\ge1$, the graph $F_2(S_n)$ is isomorphic to the total subdivision of $K_n$, where the vertices of $K_n$ are the centered configurations $\{0,i\}$, and the edge $ij$ is subdivided by the center-free configuration $\{i,j\}$.
\end{corolario}

Since subdividing edges changes neither planarity, nor the crossing number, nor the genus, Corollary~\ref{r:3.4} transfers these three parameters from $K_n$ to $F_2(S_n)$; we record the consequences in Appendix~\ref{sec:A1}.

In Appendix~\ref{sec:A} we treat several other parameters of $F_k(S_n)$, namely the chromatic number and index, the connectivity, the independence, matching and vertex-cover numbers, the planarity, the automorphism group and the spectrum. All of them are either known, as special cases of results on the token graphs of wider classes of graphs or on the doubled Johnson graphs, or almost immediate, and the point of the appendix is that the identification with $B(n;k-1,k)$ yields short and simple proofs, in contrast with the original ones, which deal with wider families of graphs.

\subsection{Distances, Wiener index, Hamiltonicity and spanning trees}\label{sec:3.3}

From here on we assume $2\le k\le(n+1)/2$, which by~\eqref{eq:1} is no loss of generality. In this subsection we compute the distances in $F_k(S_n)$ and, from them, the diameter, the radius and the Wiener index; we also determine the number of spanning trees, and we record the identification of $F_k(S_{2k-1})$ with the middle levels graph, from which its Hamiltonicity follows. We begin with the distances, which are used again in Section~\ref{sec:5}.  

\begin{lema}\label{r:3.6}
Let $X,X'$ be centered configurations, identified with their $(k-1)$-subsets of $[n]$, and let $A,A'$ be center-free ones, identified with their $k$-subsets of $[n]$. Then
\[
d(X,X')=2|X\setminus X'|,\qquad d(A,A')=2|A\setminus A'|,\qquad d(X,A)=2|X\setminus A|+1.
\]
\end{lema}

\begin{proof}
Recall that $X,X',A,A'$ are subsets of $[n]$, and they record the leaves occupied by the tokens, so the center is not among them. Every edge slides a token between the center and a leaf, so a token that moves from one leaf to another must pass through the center, and each such exchange takes two moves; moreover, every move changes the set of occupied leaves by at most one element, so the numbers of moves below cannot be improved.

To pass from $X$ to $X'$ the tokens on the leaves of $X\setminus X'$ must slide to the leaves of $X'\setminus X$, which is $|X\setminus X'|$ exchanges and hence $2|X\setminus X'|$ moves; the center-free case is identical. For a mixed pair, $\{0\}\cup X$ carries a token on the center and $A$ does not, so one further move is needed to place it on a leaf. The tokens on $X\setminus A$ slide to the leaves of $A\setminus X$ that are already matched, at two moves each, and one last move sends the token from the center to the remaining leaf of $A$. This gives $2|X\setminus A|+1$, which is $1$ when $X\subseteq A$.
\end{proof}

The eccentricities, and thus the diameter and the radius, follow from these three formulas.

\begin{corolario}\label{r:3.7}
For $2\le k\le(n+1)/2$,
\[
\operatorname{diam}\bigl(F_k(S_n)\bigr)=\min\{2k,\,n\},\qquad \operatorname{rad}\bigl(F_k(S_n)\bigr)=2k-1,
\]
and the center of $F_k(S_n)$ is the set of centered configurations when $2k\le n$, and the whole vertex set when $2k=n+1$.
\end{corolario}

\begin{proof}
We maximize each of the three formulas of Lemma~\ref{r:3.6}. Since $2k\le n+1$ we have $2(k-1)<n$ and $(k-1)+k\le n$, so two $(k-1)$-subsets of $[n]$ can always be taken disjoint, and so can a $(k-1)$-subset and a $k$-subset; hence the maxima of $|X\setminus X'|$ and of $|X\setminus A|$ are both $k-1$. For two $k$-subsets $A$ and $A'$, the maximum of $|A\setminus A'|$ is $\min\{k,n-k\}$, which equals $k$ when $2k\le n$, where disjoint $k$-subsets exist, and $k-1$ at the boundary $2k=n+1$, where they do not. Hence
\begin{align*}
\operatorname{diam}(F_k(S_n))&=\max\{\,2(k{-}1),\ 2\min\{k,n{-}k\},\ 2(k{-}1){+}1\,\}\\
&=\max\{\,2\min\{k,n{-}k\},\ 2k{-}1\,\}=\min\{2k,\,n\}.
\end{align*}
Since $\Sym(n)$ acts transitively on each side of the bipartition, all centered configurations have the same eccentricity and so do all center-free ones; by the same formulas these eccentricities are
\[
\max\{2(k-1),2(k-1)+1\}=2k-1\quad\text{and}\quad\max\{2\min\{k,n-k\},2k-1\},
\]
respectively. The second is at least the first, so the radius is $2k-1$. The second equals $2k$ when $2k\le n$, so that only the centered configurations are central, and $2k-1$ when $2k=n+1$, where consequently every vertex is central (as it must be, since in that case the complementation $\mathfrak{c}$ is an automorphism interchanging the two sides).
\end{proof}

From the formulas of Lemma~\ref{r:3.6} we can also compute the Wiener index, by a double counting argument.

\begin{teorema}\label{r:3.10}
For $1\le k\le n$, $W\bigl(F_k(S_n)\bigr)=n\binom nk\binom{n}{k-1}$.
\end{teorema}

\begin{proof}
Write $a=\binom{n}{k-1}$ and $b=\binom nk$, identify the centered configurations with their $(k-1)$-subsets $X,X'$ of $[n]$ and the center-free ones with their $k$-subsets $A,A'$. By Lemma~\ref{r:3.6},
\[
W\bigl(F_k(S_n)\bigr)\;=\;\underbrace{\sum_{\{X,X'\}}2|X\setminus X'|}_{P_1}\;+\;\underbrace{\sum_{\{A,A'\}}2|A\setminus A'|}_{P_2}\;+\;\underbrace{\sum_{X,A}\bigl(2|X\setminus A|+1\bigr)}_{P_3},
\]
where the first two sums run over the unordered pairs of distinct vertices of the same type, and the third over all $ab$ pairs of different type. The three sums are evaluated with one identity. If $\cs{F}$ and $\cs{G}$ are families of subsets of $[n]$ and $d_{\cs{F}}(x)$ denotes the number of members of $\cs{F}$ containing $x$, then counting the triples $(F,G,x)$ with $x\in F\cap G$ in two ways gives
\begin{equation}\label{eq:3}
\sum_{F\in\cs{F}}\ \sum_{G\in\cs{G}}|F\cap G|\;=\;\sum_{x\in[n]}d_{\cs{F}}(x)\,d_{\cs{G}}(x),
\end{equation}
and here every $x\in[n]$ lies in $\binom{n-1}{k-2}=\frac{k-1}{n}\,a$ of the $(k-1)$-subsets and in $\binom{n-1}{k-1}=\frac kn\,b$ of the $k$-subsets.

For $P_1$, the two ordered pairs obtained from $\{X,X'\}$ contribute equally and the $a$ diagonal terms contribute $0$, so $P_1=\sum_{(X,X')}|X\setminus X'|$ over all $a^2$ ordered pairs. Since $|X\setminus X'|=(k-1)-|X\cap X'|$ and, by~\eqref{eq:3}, $\sum_{(X,X')}|X\cap X'|=n\bigl(\frac{k-1}{n}a\bigr)^2$, and likewise for $P_2$ with $k$ and $b$ in place of $k-1$ and $a$,
\[
P_1=(k-1)a^2-\frac{(k-1)^2}{n}a^2=\frac{(k-1)(n-k+1)}{n}\,a^2,\qquad P_2=k\,b^2-\frac{k^2}{n}b^2=\frac{k(n-k)}{n}\,b^2.
\]
For $P_3$ each pair of different type occurs exactly once, the constant term contributes $ab$, and~\eqref{eq:3} gives $\sum_{X,A}|X\cap A|=n\cdot\frac{k-1}{n}a\cdot\frac knb$, so
\[
P_3=2(k-1)ab-\frac{2k(k-1)}{n}ab+ab=\frac{2(k-1)(n-k)+n}{n}\,ab.
\]
It remains to add the three. Using $b=\frac{n-k+1}{k}\,a$ to express everything in terms of $a$, and writing $m=n-k$, we obtain
\[
P_1+P_2+P_3=\frac{(n-k+1)\,a^2}{nk}\Bigl[\,k(k-1)+m(m+1)+2(k-1)m+n\,\Bigr],
\]
and the bracket equals $(k+m)^2-(k+m)+n=n^2$, since $k+m=n$. Therefore \[W\bigl(F_k(S_n)\bigr)=n\,a\cdot\frac{(n-k+1)a}{k}=n\,a\,b=n\binom nk\binom{n}{k-1}.\]
\end{proof}

For Hamiltonicity we relate $F_k(S_n)$ to a graph that has been studied in its own right. For $m\ge1$, the \emph{middle levels graph} $M_m$ has as vertices the subsets of $[2m+1]$ of size $m$ and those of size $m+1$ (corresponding to the two middle levels of the subset lattice of $[2m+1]$) with $S$ and $T$ adjacent whenever $S\subset T$. In the notation of this subsection, $M_m\cong B(2m+1;\,m,\,m+1)$; it is also the \emph{doubled odd graph} $2O_{m+1}$, the bipartite double of the odd graph $O_{m+1}$ (see~\cite[Section~9.1D]{b11}), under the correspondence that sends an $(m+1)$-set to its complement. By a theorem of M\"utze~\cite{b43} (see also the shorter proof of Gregor, M\"utze and Nummenpalo~\cite{b30}), $M_m$ is Hamiltonian for every $m\ge1$. Comparing the two descriptions and given that $2k-1=2(k-1)+1$, Theorem~\ref{r:3.1} gives
\begin{equation}\label{eq:2}
F_k(S_{2k-1})\;\cong\;B(2k-1;k-1,k)\;\cong\;M_{k-1}\;\cong\;2O_k.
\end{equation}
The same object thus has three names, and all three will be used, namely $F_k(S_{2k-1})$ in Sections~\ref{sec:5} and~\ref{sec:6}, the middle levels graph here, and the doubled odd graph in Subsection~\ref{sec:6.2}, where perfect codes in $O_k$ enter.

\begin{teorema}\label{r:3.9}
For $2\le k\le(n+1)/2$, the graph $F_k(S_n)$ is Hamiltonian if and only if $2k=n+1$.
\end{teorema}

\begin{proof}
In order to have a Hamiltonian cycle, a bipartite graph must have both sides of the same order. Applied to $F_k(S_n)$ this requires $\binom{n}{k-1}=\binom nk$, that is $k-1=n-k$, or $2k=n+1$. Conversely, if $n=2k-1$ then $F_k(S_n)\cong M_{k-1}$ by~\eqref{eq:2}, and $M_{k-1}$ is Hamiltonian~\cite{b43}.
\end{proof}

Finally we turn to the Laplacian spectrum. Recall that the Laplacian matrix of a graph $G$ is $L(G)=D(G)-A(G)$, where $A(G)$ is the adjacency matrix and $D(G)$ the diagonal matrix of degrees; its eigenvalues are nonnegative, the smallest one is $0$, and the second smallest is the algebraic connectivity $a(G)$. Dalf\'o et al.~\cite[Section~5]{b18} determined the Laplacian eigenvalues of $F_k(S_n)$ from the quotient matrix of its bipartition, with their multiplicities when $n=2k-1$, and proved that $a(F_k(S_n))=1$~\cite[Theorem~7.2]{b18}. The next proposition gives the multiplicities for all $n$ and $k$. We prove it in Appendix~\ref{sec:A5}, together with the adjacency spectrum, from the singular values of the inclusion matrix of $(k-1)$-sets against $k$-sets. From the Laplacian spectrum, Kirchhoff's matrix-tree theorem gives the number of spanning trees.

\begin{proposicion}[see~\cite{b18}]\label{r:3.13}
For $2\le k\le(n+1)/2$, the Laplacian spectrum of $F_k(S_n)$ consists of the integers
\[
j\quad\text{and}\quad n+1-j,\qquad\text{each with multiplicity}\quad\binom nj-\binom{n}{j-1},\qquad 0\le j\le k-1,
\]
together with $k$ with multiplicity $\binom nk-\binom{n}{k-1}$. Consequently $a\bigl(F_k(S_n)\bigr)=1$ for $2\le k\le n-1$.
\end{proposicion}

\begin{corolario}\label{r:sptrees}
For $2\le k\le(n+1)/2$ the number of spanning trees of $F_k(S_n)$ is
\[
\tau_{\mathrm{sp}}\bigl(F_k(S_n)\bigr)=\frac{n+1}{\binom{n+1}{k}}\;k^{\binom nk-\binom{n}{k-1}}\;\prod_{j=1}^{k-1}\bigl(j(n+1-j)\bigr)^{\binom nj-\binom{n}{j-1}}.
\]
\end{corolario}

\begin{proof}
By Proposition~\ref{r:3.13}, the Laplacian eigenvalues of $F_k(S_n)$ are $j$ and $n+1-j$, each with multiplicity $\binom nj-\binom{n}{j-1}$, for $0\le j\le k-1$, together with $k$ with multiplicity $\binom nk-\binom{n}{k-1}$. Kirchhoff's matrix-tree theorem~\cite[Proposition~1.3.4]{bBH} gives $\tau_{\mathrm{sp}}=|V|^{-1}\prod_{\mu\ne0}\mu$, with $|V|=\binom{n+1}{k}$ and $\mu$ running over the nonzero Laplacian eigenvalues with multiplicity. The index $j=0$ produces the pair $0$ and $n+1$, each once, and $0$ is the eigenvalue omitted from the product, so only the factor $n+1$ survives; each index $1\le j\le k-1$ contributes $\bigl(j(n+1-j)\bigr)^{\binom nj-\binom{n}{j-1}}$, and the eigenvalue $k$ contributes $k^{\binom nk-\binom{n}{k-1}}$.
\end{proof}

\section{Domination}\label{sec:4}
The aim of this section is to show Theorem~\ref{r:1.2}, which is the union of Theorem~\ref{r:4.5} (the upper bound), Theorem~\ref{r:4.8} (the lower bound), Theorem~\ref{r:4.10} (the exact bound for $k=3$) and Theorem~\ref{r:4.12} (the upper bound for $k=4$). Subsection~\ref{sec:4.1} translates domination into an extremal question about $(k-1)$-uniform hypergraphs; the upper bound comes from a construction (Subsection~\ref{sec:4.2}) and the lower bound from the removal lemma (Subsection~\ref{sec:4.3}); Subsections~\ref{sec:4.4} and~\ref{sec:4.5} treat $k\le4$ and the comparison with~\cite{b23}.

\subsection{From dominating sets to hypergraphs}\label{sec:4.1}

For $k=2$ the domination number is known exactly, $\gamma(F_2(S_n))=n-1$~\cite{b23} (see Subsection~\ref{sec:4.4}), so throughout this section $k\ge 3$ is fixed, and we work in $B(n;k-1,k)$ via Theorem~\ref{r:3.1}. The starting point is that a dominating set is determined by what it leaves out on the side of the $(k-1)$-sets, which is the smaller side since $k\le (n+1)/2$, and that this remainder is best read as a hypergraph.

Given a candidate set $D$, we write $D_{k-1}=D\cap\binom{[n]}{k-1}$ and $D_k=D\cap\binom{[n]}{k}$, and we denote by $\mathcal{H}=\mathcal{H}(D)$ the $(k-1)$-uniform hypergraph on $[n]$ whose edges are the $(k-1)$-sets that are \emph{not} in $D_{k-1}$; as usual, $e(\mathcal{H})$ denotes the number of edges of $\mathcal{H}$. We say that a $k$-set $A$ is \emph{complete} (with respect to $\mathcal{H}$) if its $k$ subsets of size $k-1$ are all edges of $\mathcal{H}$; equivalently, if $A$ generates a copy of $K_k^{(k-1)}$, the complete $(k-1)$-uniform hypergraph on $k$ vertices. We write $K(\mathcal{H})$ for the set of complete $k$-sets.

\begin{proposicion}\label{r:4.1}
The set $D=D_{k-1}\cup D_k$ is dominating in $B(n;k-1,k)$ if and only if
\begin{enumerate}[label=(\roman*)]
\item every complete $k$-set of $\mathcal{H}$ belongs to $D_k$, and
\item every edge of $\mathcal{H}$ is contained in some member of $D_k$.
\end{enumerate}
Consequently,
\[
\gamma\bigl(F_k(S_n)\bigr)=\min_{\mathcal{H}}\left[\binom{n}{k-1}-e(\mathcal{H})+d_k(\mathcal{H})\right],
\]
where $\mathcal{H}$ ranges over all $(k-1)$-uniform hypergraphs on $[n]$ and $d_k(\mathcal{H})$ denotes the minimum size of a family of $k$-sets satisfying (i) and (ii).
\end{proposicion}

\begin{proof}
The neighborhood of a $k$-set $A$ consists exactly of its $k$ faces of size $k-1$. Hence $A$ is dominated if and only if $A\in D_k$ or one of its faces is in $D_{k-1}$. By the definition of $\mathcal{H}$, a face of $A$ fails to lie in $D_{k-1}$ exactly when it is an edge of $\mathcal{H}$. Hence $A$ has no face in $D_{k-1}$ if and only if all $k$ of its faces are edges of $\mathcal{H}$, that is, if and only if $A$ is complete. Every $k$-set that is not complete is therefore already dominated by one of its faces, and a complete one is dominated only if it belongs to $D_k$, as stated in condition (i).

The neighborhood of a $(k-1)$-set $X$ consists exactly of the $k$-sets that contain it. Hence $X$ is dominated if and only if $X\in D_{k-1}$ or $X$ is covered by a member of $D_k$; since the $(k-1)$-sets not in $D_{k-1}$ are exactly the edges of $\mathcal{H}$, all of them are dominated if and only if (ii) holds.

Note that $|D_{k-1}|=\binom{n}{k-1}-e(\mathcal{H})$, so $|D|=\binom{n}{k-1}-e(\mathcal{H})+|D_k|$. The minimum is now reached from both sides. If $D$ is a smallest dominating set and $\mathcal{H}=\mathcal{H}(D)$ is its hypergraph, then $D_k$ satisfies (i) and (ii), so $|D_k|\ge d_k(\mathcal{H})$ and $\gamma(F_k(S_n))\ge\binom{n}{k-1}-e(\mathcal{H})+d_k(\mathcal{H})$, which is at least the minimum over all $\mathcal{H}$. Conversely, given any $(k-1)$-uniform hypergraph $\mathcal{H}'$, take $D'_k$ of minimum size satisfying (i) and (ii) for $\mathcal{H}'$ and let $D'_{k-1}$ be the $(k-1)$-sets that are not edges of $\mathcal{H}'$; the hypergraph associated with $D'=D'_{k-1}\cup D'_k$ is $\mathcal{H}'$ itself, so $D'$ is dominating and $\gamma(F_k(S_n))\le\binom{n}{k-1}-e(\mathcal{H}')+d_k(\mathcal{H}')$.
\end{proof}

\begin{lema}\label{r:4.2}
Let $\mathcal{H}$ be a $(k-1)$-uniform hypergraph on $[n]$, let $m_1$ be the number of edges of $\mathcal{H}$ contained in at least one complete $k$-set, and let $m_0=e(\mathcal{H})-m_1$. Then
\[
d_k(\mathcal{H})\ \ge\ |K(\mathcal{H})|+\left\lceil\frac{m_0}{k-1}\right\rceil,
\]
and if $m_0=0$ then $d_k(\mathcal{H})=|K(\mathcal{H})|$. Moreover, $|K(\mathcal{H})|\ \ge\ m_1/k$.
\end{lema}

\begin{proof}
Condition (i) of Proposition~\ref{r:4.1} forces $K(\mathcal{H})\subseteq D_k$. By condition (ii), every edge of $\mathcal{H}$ that lies in no complete $k$-set must be covered by a non-complete member of $D_k$; and a non-complete $k$-set contains at most $k-1$ edges of $\mathcal{H}$, since all $k$ of them would make it complete. This gives the first inequality. If $m_0=0$, the family $D_k=K(\mathcal{H})$ satisfies (i) and (ii) and attains the lower bound.

For the last inequality, we count incidences between complete $k$-sets and their faces. Each complete $k$-set has exactly $k$ faces, all of them edges of $\mathcal{H}$, and each of the $m_1$ covered edges lies in at least one complete $k$-set.
\end{proof}

For $k=3$ the hypergraph $\mathcal{H}$ is an ordinary graph, its complete $3$-sets are exactly its triangles, and the two bounds of Lemma~\ref{r:4.2} read $d_3(\mathcal{H})\ge|K(\mathcal{H})|+\lceil m_0/2\rceil$ and $|K(\mathcal{H})|\ge m_1/3$.

\subsection{An upper bound: the book construction}\label{sec:4.2}

Taking $\mathcal{H}$ empty in Proposition~\ref{r:4.1} gives the trivial bound $\gamma(F_k(S_n))\le\binom{n}{k-1}$, which is attained by the set of all $(k-1)$-sets. Every $(k-1)$-uniform hypergraph $\mathcal{H}$ on $[n]$ improves on it by
\[
f(\mathcal{H})=e(\mathcal{H})-d_k(\mathcal{H}),
\]
a quantity we call the \emph{saving} of $\mathcal{H}$, since each edge removes one $(k-1)$-set from the dominating set, and $d_k(\mathcal{H})$ is the price paid in $k$-sets to keep it dominating. Bounding $\gamma(F_k(S_n))$ from above is therefore the problem of making the saving large, that is, of exhibiting a hypergraph with many edges and few complete $k$-sets. We now give a construction with these properties.

\begin{construccion}\label{r:4.3}
Let $k\ge 3$ and let $s,p$ be integers with $s\ge 1$, $p\ge k-2$ and $2s+p=n$. Partition
\[
[n]=\{a_1,b_1,\dots,a_s,b_s\}\cup P,\qquad |P|=p,
\]
and call the pairs $\{a_i,b_i\}$ the \emph{spines} and the elements of $P$ the \emph{pages}. Write $S=\{a_1,b_1,\dots,a_s,b_s\}$ for the set of spine vertices, so that $[n]=S\cup P$. Let $\mathcal{H}=\mathcal{H}_{s,p}$ be the $(k-1)$-uniform hypergraph on $[n]$ whose edges are
\begin{enumerate}[label=(\Roman*)]
\item the sets $\{a_i,b_i\}\cup Q$ with $1\le i\le s$ and $Q\in\binom{P}{k-3}$, and
\item the sets $\{u\}\cup Q$ with $u\in\{a_i,b_i\}$ for some $1\le i\le s$ and $Q\in\binom{P}{k-2}$.
\end{enumerate}
Counting by type, $e(\mathcal{H})=s\binom{p}{k-3}+2s\binom{p}{k-2}$.
\end{construccion}

\begin{figure}[ht]
\centering
\includegraphics[width=0.85\textwidth]{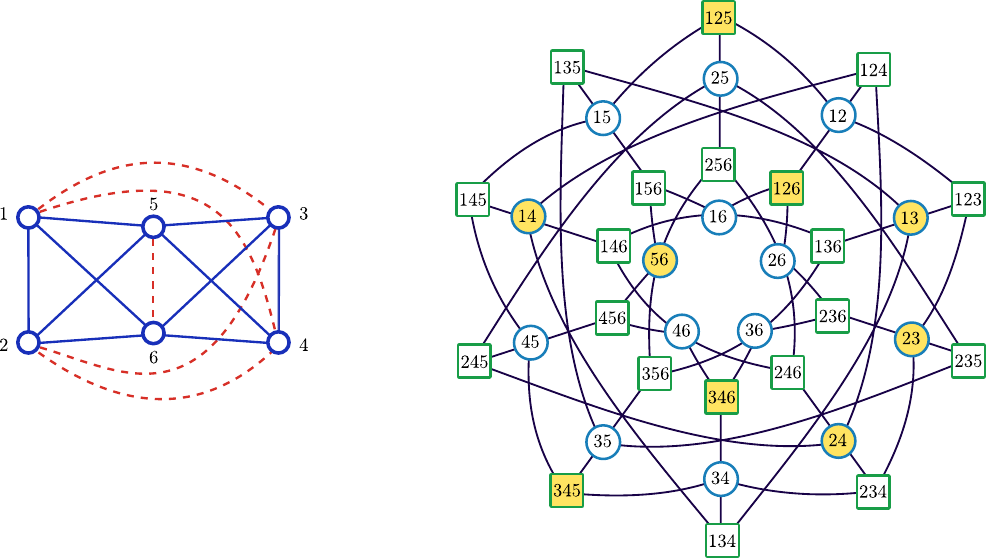}
\caption{The book construction for $F_3(S_6)$, with spines $\{1,2\}$ and $\{3,4\}$ and pages $\{5,6\}$. On the left, the graph $\mathcal{H}$, its non-edges drawn dashed in red. On the right, the resulting dominating set of $F_3(S_6)$, in yellow, consisting of the five non-edges of $\mathcal{H}$ together with its four triangles.}
\label{fig:2}
\end{figure}

The name comes from the case $k=3$, where $\mathcal{H}$ is an ordinary graph whose edges are the $s$ spines and the $2sp$ spine--page pairs, so that $\mathcal{H}$ is a union of $s$ \emph{books}, each an edge together with a set of vertices joined to both of its ends. The $s$ books have disjoint spines and share the same pages. For larger $k$ the pages enter $k-2$ at a time, but the shape is unchanged. See an example in Figure~\ref{fig:2}.

The next lemma describes the complete $k$-sets of $\mathcal{H}_{s,p}$ and shows that every edge lies in one of them.

\begin{lema}\label{r:4.4}
Let $\mathcal{H}=\mathcal{H}_{s,p}$ be as in Construction~\ref{r:4.3}. Then:
\begin{enumerate}[label=(\roman*)]
\item the complete $k$-sets of $\mathcal{H}$ are the sets $\{a_i,b_i\}\cup Q$ with $1\le i\le s$ and $Q\in\binom{P}{k-2}$, so $|K(\mathcal{H})|=s\binom{p}{k-2}$;
\item every edge of $\mathcal{H}$ is contained in a complete $k$-set, that is, $m_0=0$.
\end{enumerate}
\end{lema}

\begin{proof}
Let $A=\{a_i,b_i\}\cup Q$ with $|Q|=k-2$, so $|A|=k$. Its faces are $\{a_i\}\cup Q$ and $\{b_i\}\cup Q$, of type (II), and the sets $\{a_i,b_i\}\cup(Q\setminus\{q\})$ with $q\in Q$, of type (I); all are edges, so $A$ is complete.

Conversely, let $A$ be complete and $t=|A\cap S|$. From Construction~\ref{r:4.3} we obtain that every edge meets $S$ in one or two vertices, and in two only when these are the two ends of a same spine. If $t\le 1$, then $A$ has a face contained in $P$, but such a face is not an edge. Suppose $t\ge 3$. If $A\cap P\ne\varnothing$, removing a page leaves a face with $t\ge 3$ vertices of $S$, again not an edge. If $A\cap P=\varnothing$, then every face of $A$ has $k-1\ge 3$ vertices of $S$ when $k\ge 4$, and is not an edge; when $k=3$ the three faces of $A$ would have to be spines, which is impossible since each vertex lies in exactly one spine. Hence $t=2$.

Write $A\cap S=\{u,v\}$, so that $|A\cap P|=k-2\ge 1$. Removing a page from $A$ leaves a face containing both $u$ and $v$, which is an edge only if $\{u,v\}$ is a spine. Therefore $A=\{a_i,b_i\}\cup Q$ with $Q=A\cap P$ of size $k-2$. Moreover, distinct pairs $(i,Q)$ give distinct sets, so $|K(\mathcal{H})|=s\binom{p}{k-2}$. This proves (i).

We now show (ii). An edge of type (I), say $\{a_i,b_i\}\cup Q$ with $|Q|=k-3$, is contained in $\{a_i,b_i\}\cup(Q\cup\{y\})$ for any page $y\notin Q$, and such a page exists because $p\ge k-2>k-3=|Q|$. An edge of type (II), say $\{u\}\cup W$ with $u\in\{a_i,b_i\}$ and $|W|=k-2$, is contained in $\{a_i,b_i\}\cup W$.
\end{proof}

With Lemma~\ref{r:4.4} in hand, the saving of the books can be computed exactly, and choosing the number of spines optimally gives the general upper bound of Theorem~\ref{r:1.2}.

\begin{teorema}\label{r:4.5}
With $s$ and $p$ as in Construction~\ref{r:4.3},
\[
\gamma\bigl(F_k(S_n)\bigr)\ \le\ \binom{n}{k-1}-s\binom{p+1}{k-2}.
\]
Optimizing over $s$ and $p$ one obtains, as $n\to\infty$,
\[
\gamma\bigl(F_k(S_n)\bigr)\ \le\ \left(1-\frac12\Bigl(\frac{k-2}{k-1}\Bigr)^{k-2}+o(1)\right)\binom{n}{k-1}.
\]
\end{teorema}

\begin{proof}
By Lemma~\ref{r:4.4} the hypergraph $\mathcal{H}=\mathcal{H}_{s,p}$ of Construction~\ref{r:4.3} has $m_0=0$ and $|K(\mathcal{H})|=s\binom{p}{k-2}$, so Lemma~\ref{r:4.2} gives $d_k(\mathcal{H})=|K(\mathcal{H})|$. By Proposition~\ref{r:4.1} and using Pascal's identity $\binom{p}{k-3}+\binom{p}{k-2}=\binom{p+1}{k-2}$,
\[
\gamma\bigl(F_k(S_n)\bigr)\le\binom{n}{k-1}-e(\mathcal{H})+|K(\mathcal{H})|=\binom{n}{k-1}-s\left[\binom{p}{k-3}+\binom{p}{k-2}\right]=\binom{n}{k-1}-s\binom{p+1}{k-2}.
\]

It remains to optimize. Here and below we use the elementary estimate that, for a fixed $r\ge 0$ and $m\to\infty$,
\begin{equation}\label{eq:4}
\binom{m}{r}=\frac{m(m-1)\cdots(m-r+1)}{r!}=\frac{m^r}{r!}\prod_{i=1}^{r-1}\Bigl(1-\frac{i}{m}\Bigr)=\bigl(1+o(1)\bigr)\frac{m^r}{r!},
\end{equation}
since the product has a fixed number of factors, each tending to $1$. Write $p=\alpha n$ with $\alpha\in(0,1)$, so $s=(n-p)/2=(1-\alpha)n/2$. Applying \eqref{eq:4} with $m=p+1$ and $r=k-2$, and using $p+1=\bigl(1+o(1)\bigr)\alpha n$,
\[
s\binom{p+1}{k-2}=\bigl(1+o(1)\bigr)\,\frac{(1-\alpha)\,\alpha^{k-2}}{2\,(k-2)!}\,n^{k-1}.
\]
The function $\varphi(\alpha)=(1-\alpha)\alpha^{k-2}$ has derivative $\varphi'(\alpha)=\alpha^{k-3}\bigl[(k-2)-(k-1)\alpha\bigr]$, which vanishes at $\alpha^*=\frac{k-2}{k-1}$ and changes sign from positive to negative there, so this is the maximum, where $\varphi(\alpha^*)=\frac{1}{k-1}\bigl(\frac{k-2}{k-1}\bigr)^{k-2}$. Dividing by $\binom{n}{k-1}=\bigl(1+o(1)\bigr)n^{k-1}/(k-1)!$, which is \eqref{eq:4} with $m=n$ and $r=k-1$, we obtain
\[
\frac{s\binom{p+1}{k-2}}{\binom{n}{k-1}}\ \longrightarrow\ \frac{(k-1)!}{2\,(k-2)!\,(k-1)}\Bigl(\frac{k-2}{k-1}\Bigr)^{k-2}\ =\ \frac12\Bigl(\frac{k-2}{k-1}\Bigr)^{k-2},
\]
which gives the asymptotic statement.
\end{proof}

The spine of Construction~\ref{r:4.3} could be taken of any size $\ell\ge2$, the edges being the $(k-1)$-subsets of the sets $S_i\cup Q$ with $S_i$ a spine and $Q\in\binom{P}{k-\ell}$. For $\ell\ge3$, however, such a hypergraph has at most $k\,s\binom{p}{k-\ell}=O(n^{k-\ell+1})=o(n^{k-1})$ edges, so its saving $f(\mathcal{H})\le e(\mathcal{H})$ is $o\bigl(\binom{n}{k-1}\bigr)$ regardless of the choice of $D_k$. Thus, only the edge spine yields a saving of order $\binom{n}{k-1}$.

\subsection{A lower bound: the removal lemma}\label{sec:4.3}

For an $r$-uniform hypergraph $\mathcal{F}$, let $\ex(n;\mathcal{F})$ be the maximum number of edges of an $\mathcal{F}$-free $r$-uniform hypergraph on $n$ vertices, and let $\pi(\mathcal{F})=\lim_{n\to\infty}\ex(n;\mathcal{F})/\binom{n}{r}$ be its Turán density. The hypergraph that matters here is $K_k^{(k-1)}$, with the same $k$ as the number of tokens, and we abbreviate $\pi_k:=\pi\bigl(K_k^{(k-1)}\bigr)$. 
By Mantel's theorem~\cite{bMantel}, $\pi_3=\pi(K_3)=\frac12$. For $k\ge 4$ the value of $\pi_k$ is not known. Turán~\cite{b54} conjectured $\pi_4=\frac59$, and the best known upper bound is $\pi_4\le 0.5578$, due to Jeong, Park, Im, Lee and Yang~\cite{b35}, who improved earlier flag-algebra bounds of Razborov~\cite{b48} and Baber~\cite{bBaber}. For general $k$ we use the bound of de Caen~\cite{b20},
\begin{equation}\label{eq:5}
\pi\bigl(K_{r+1}^{(r)}\bigr)\ \le\ 1-\frac1r\qquad(r\ge 2).
\end{equation}

\begin{teorema}[Removal lemma~\cite{bRSz,bNRS,bRS06,bGowers}]\label{r:rem}
Let $\mathcal{F}$ be an $r$-uniform hypergraph with $v(\mathcal{F})$ vertices. For every $\varepsilon>0$ there exist $\delta>0$ and $n_0$ such that every $r$-uniform hypergraph on $n\ge n_0$ vertices with at most $\delta\,n^{v(\mathcal{F})}$ copies of $\mathcal{F}$ can be made $\mathcal{F}$-free by removing at most $\varepsilon\,n^{r}$ edges.
\end{teorema}

For $r=2$ and $\mathcal{F}=K_3$ this is the triangle removal lemma of Ruzsa and Szemerédi~\cite{bRSz}. For $\mathcal{F}=K_{r+1}^{(r)}$, which is the case we use, it was proved by Nagle, Rödl and Schacht~\cite{bNRS}; Rödl and Skokan~\cite{bRS06} extended it to every $\mathcal{F}$, and Gowers~\cite{bGowers} obtained it independently for every $\mathcal{F}$.

Here, a copy of $\mathcal{F}=K_k^{(k-1)}$ in a $(k-1)$-uniform hypergraph $\mathcal{H}$ is a $k$-set all of whose $(k-1)$-subsets are edges, that is, a complete $k$-set; thus the number of copies of $K_k^{(k-1)}$ in $\mathcal{H}$ is $|K(\mathcal{H})|$.

\begin{teorema}\label{r:4.8}
Let $k\ge 3$ be fixed. Then
\[
\gamma\bigl(F_k(S_n)\bigr)\ \ge\ \Bigl(1-\frac{k-2}{k-1}\,\pi_k-o(1)\Bigr)\binom{n}{k-1}.
\]
\end{teorema}

\begin{proof}
Fix $\varepsilon>0$ and let $\delta$ and $n_0$ be given by Theorem~\ref{r:rem} for $r=k-1$ and $\mathcal{F}=K_k^{(k-1)}$, which has $v(\mathcal{F})=k$ vertices. By the definition of $\pi_k$ there is $n_1$ such that $\ex(n;K_k^{(k-1)})\le(\pi_k+\varepsilon)\binom{n}{k-1}$ for $n\ge n_1$. Let $n\ge\max\{n_0,n_1\}$ be large enough that, in addition, $\delta n^k>\binom{n}{k-1}$.

Let $D=D_{k-1}\cup D_k$ be a dominating set of minimum size of $B(n;k-1,k)$ and let $\mathcal{H}=\mathcal{H}(D)$, so that $|D|=\binom{n}{k-1}-e(\mathcal{H})+|D_k|$ and $D_k$ satisfies conditions (i) and (ii) of Proposition~\ref{r:4.1}. If $|K(\mathcal{H})|>\delta n^k$, then by (i) of Proposition~\ref{r:4.1},
\[
|D|\ \ge\ |D_k|\ \ge\ |K(\mathcal{H})|\ >\ \delta n^k\ >\ \binom{n}{k-1},
\]
which contradicts the trivial bound $\gamma(F_k(S_n))\le\binom{n}{k-1}$. Hence $|K(\mathcal{H})|\le\delta n^k$, and by Theorem~\ref{r:rem} there is $R\subseteq E(\mathcal{H})$ with $|R|\le\varepsilon n^{k-1}$ such that $\mathcal{H}'=\mathcal{H}-R$ has no complete $k$-set. This has two consequences.
\begin{enumerate}[label=(\alph*)]
\item $\mathcal{H}'$ is $K_k^{(k-1)}$-free, so $e(\mathcal{H}')\le(\pi_k+\varepsilon)\binom{n}{k-1}$, and hence $e(\mathcal{H})=e(\mathcal{H}')+|R|\le(\pi_k+\varepsilon)\binom{n}{k-1}+\varepsilon n^{k-1}$.
\item By (ii), every edge of $\mathcal{H}'$ is contained in some member of $D_k$. Since $\mathcal{H}'$ has no complete $k$-set, every $k$-set contains at most $k-1$ edges of $\mathcal{H}'$, so that $(k-1)|D_k|\ge e(\mathcal{H}')\ge e(\mathcal{H})-\varepsilon n^{k-1}$.
\end{enumerate}
By (b) and then by (a),
\begin{align*}
|D|\ &\ge\ \binom{n}{k-1}-e(\mathcal{H})+\frac{e(\mathcal{H})-\varepsilon n^{k-1}}{k-1}\ =\ \binom{n}{k-1}-\frac{k-2}{k-1}\,e(\mathcal{H})-\frac{\varepsilon n^{k-1}}{k-1}\\
&\ge\ \Bigl(1-\frac{k-2}{k-1}\,\pi_k-\frac{k-2}{k-1}\,\varepsilon\Bigr)\binom{n}{k-1}-\varepsilon n^{k-1}.
\end{align*}
Since $n^{k-1}\le(k-1)^{k-1}\binom{n}{k-1}$, which follows from $\binom{m}{r}=\prod_{i=0}^{r-1}\frac{m-i}{r-i}\ge\bigl(\frac mr\bigr)^r$, we obtain
\[
\gamma\bigl(F_k(S_n)\bigr)\ =\ |D|\ \ge\ \Bigl(1-\frac{k-2}{k-1}\,\pi_k-\bigl(1+(k-1)^{k-1}\bigr)\varepsilon\Bigr)\binom{n}{k-1}
\]
for every sufficiently large $n$. Since $\varepsilon>0$ is arbitrary, the theorem follows.
\end{proof}

For $k=3$ the argument reads as follows. There $\mathcal{H}'$ is a triangle-free graph, so $e(\mathcal{H}')\le n^2/4$ by Mantel's theorem~\cite{bMantel}, and every triple of $D_3$ contains at most two edges of $\mathcal{H}'$; hence
\[
|D|\ \ge\ \binom n2-e(\mathcal{H})+\frac{e(\mathcal{H}')}{2}\ =\ \binom n2-\frac{e(\mathcal{H}')}{2}-|R|\ \ge\ \binom n2-\frac{n^2}{8}-o(n^2).
\]

Applying de Caen's bound \eqref{eq:5} with $r=k-1$ we obtain $\pi_k\le\frac{k-2}{k-1}$, so $1-\frac{k-2}{k-1}\pi_k\ge1-\bigl(\frac{k-2}{k-1}\bigr)^2=\frac{2k-3}{(k-1)^2}$, and by Theorem~\ref{r:4.8} we obtain:

\begin{corolario}\label{r:4.9}
Let $k\ge 3$ be fixed. Then
\[
\gamma\bigl(F_k(S_n)\bigr)\ \ge\ \Bigl(\frac{2k-3}{(k-1)^2}-o(1)\Bigr)\binom{n}{k-1}.
\]
\end{corolario}

In particular, Mantel's theorem~\cite{bMantel} gives $\pi_3=\frac12$ and the bound $\pi_4\le 0.5578$ of~\cite{b35} gives
\begin{equation}\label{eq:6}
\gamma\bigl(F_3(S_n)\bigr)\ \ge\ \Bigl(\tfrac34-o(1)\Bigr)\binom{n}{2},\qquad\text{and}\qquad\gamma\bigl(F_4(S_n)\bigr)\ \ge\ \bigl(0.628-o(1)\bigr)\binom{n}{3}.
\end{equation}

\subsection{The cases \texorpdfstring{$k\le 4$}{k <= 4}}\label{sec:4.4}

For $k=2$ the problem is settled. Fabila-Monroy and Gómez-Galicia~\cite{b23} showed that
\[
\gamma\bigl(F_2(S_n)\bigr)=n-1\qquad(n\ge 2),
\]
and Theorem~\ref{r:3.1} makes the upper bound evident, since by Corollary~\ref{r:3.4} the graph is the total subdivision of $K_n$, so a minimum dominating set is formed by $n-2$ branch vertices, together with the subdivision vertex of the pair formed by the two omitted ones. For $k=3$ the general bounds above can be made explicit, and for $k=4$ a different construction improves the upper bound.

For $k=3$, Theorem~\ref{r:4.5} reads $\gamma(F_3(S_n))\le\binom{n}{2}-s(p+1)$, and choosing $s$ optimally gives a closed form.

\begin{teorema}\label{r:4.10}
For $n\ge 3$,
\[
\gamma\bigl(F_3(S_n)\bigr)\ \le\ \binom{n}{2}-\max_{\substack{s,p\ge 1\\ 2s+p=n}}s(p+1)\ =\ \binom{n}{2}-\left\lfloor\frac{(n+1)^2}{8}\right\rfloor.
\]
\end{teorema}

\begin{proof}
The inequality is Theorem~\ref{r:4.5} for $k=3$, so only the maximization is left. Note that it has to be carried out over the integers, since Theorem~\ref{r:4.5} optimizes over the reals and returns only the asymptotic proportion of pages, here $p\sim n/2$, whereas the closed form above is exact for every $n$.

Writing $p=n-2s$ makes the saving a function of $s$ alone, namely $g(s)=s(n-2s+1)$, to be maximized over the integers $1\le s\le(n-1)/2$, the range being forced by $p\ge 1$. Completing the square,
\[
g(s)=\frac{(n+1)^2}{8}-2\left(s-\frac{n+1}{4}\right)^2,
\]
so $g$ is maximized at an integer $s^*$ closest to $(n+1)/4$, and $g(s^*)=\frac{1}{8}(n+1)^2-2d^2$, where $d\le\frac12$ is the distance from $(n+1)/4$ to $s^*$. Since $g(s^*)$ is an integer and $0\le 2d^2\le\frac12<1$, it follows that $g(s^*)=\lfloor(n+1)^2/8\rfloor$. Finally, $s^*\ge1$ and $p=n-2s^*\ge1$ for every $n\ge3$, so $s^*$ is admissible.
\end{proof}

\begin{corolario}\label{r:4.11}
As $n\to\infty$,
\[
\Bigl(\tfrac34-o(1)\Bigr)\binom{n}{2}\ \le\ \gamma\bigl(F_3(S_n)\bigr)\ \le\ \binom{n}{2}-\left\lfloor\frac{(n+1)^2}{8}\right\rfloor\ =\ \Bigl(\tfrac34+o(1)\Bigr)\binom{n}{2};
\]
in particular, $\gamma(F_3(S_n))=\bigl(\frac34+o(1)\bigr)\binom n2$.
\end{corolario}

Integer programming shows that the bound of Theorem~\ref{r:4.10} is attained for all $3\le n\le 10$, where $\gamma(F_3(S_n))$ takes the values $1,3,6,9,13,18,24,30$. By Corollary~\ref{r:4.11} the book construction is asymptotically optimal; it is exact for all these values, and we conjecture that this persists (Conjecture~\ref{r:7.2} in Section~\ref{sec:7}).

For $k=4$, Theorem~\ref{r:4.5} reads
\[
\gamma\bigl(F_4(S_n)\bigr)\ \le\ \binom{n}{3}-s\,\frac{p(p+1)}{2}\ =\ \Bigl(\tfrac79+o(1)\Bigr)\binom{n}{3},
\]
which already improves the bound $\bigl(\frac56+o(1)\bigr)\binom{n}{3}$ of~\cite{b23} (the comparison with~\cite{b23} for general $k$ is made in Subsection~\ref{sec:4.5}), but the book construction is not optimal. Proposition~\ref{r:4.1} asks for a $3$-uniform hypergraph $\mathcal{H}$ and a family $D_4$ of $4$-sets that covers the edges of $\mathcal{H}$ and contains its complete $4$-sets. If $\mathcal{H}$ has no complete $4$-set and its edges are partitioned by $4$-sets that contain exactly three of them each, then $|D_4|=e(\mathcal{H})/3$ and $\gamma(F_4(S_n))\le\binom n3-\frac23\,e(\mathcal{H})$. Turán's construction~\cite{b54}, of density $\frac59$, admits such a partition up to $O(n^2)$ edges.

\begin{teorema}\label{r:4.12}
As $n\to\infty$, $\gamma\bigl(F_4(S_n)\bigr)\le\bigl(\frac{17}{27}+o(1)\bigr)\binom{n}{3}$.
\end{teorema}
\begin{proof}
Assume first that $n=3q$ for some positive integer $q$, and partition $[n]=V_0\cup V_1\cup V_2$ into three parts of size $q$; indices are taken modulo $3$. Let $\mathcal{H}$ be the $3$-uniform hypergraph whose edges are the triples $\{a,b,c\}$ with $a\in V_0$, $b\in V_1$ and $c\in V_2$, and the triples $\{x,x',y\}$ with $x,x'\in V_i$ and $y\in V_{i+1}$ for some $i\in\{0,1,2\}$. Then
\[
e(\mathcal{H})=q^3+3\binom{q}{2}q=\Bigl(\frac{1}{27}+\frac{1}{18}+o(1)\Bigr)n^3=\Bigl(\frac59+o(1)\Bigr)\binom{n}{3}.
\]
The hypergraph $\mathcal{H}$ has no complete $4$-set. Indeed, a $4$-set $Q$ has two points $x,x'$ in some $V_i$; the two faces of $Q$ that contain $x$ and $x'$ are edges only if the other two points $y,y'$ of $Q$ lie in $V_{i+1}$, but then the face $\{x,y,y'\}$, with two points in $V_{i+1}$ and one in $V_i$, is not an edge. Hence condition~(i) of Proposition~\ref{r:4.1} holds for every family of $4$-sets, and it suffices to cover the edges of $\mathcal{H}$ with few $4$-sets.

We use $4$-sets of two kinds:
\begin{itemize}
    \item a set $\{a,a',b,c\}$ with $a,a'\in V_0$, $b\in V_1$ and $c\in V_2$ contains the three edges $\{a,a',b\}$, $\{a,b,c\}$ and $\{a',b,c\}$;
    \item a set $\{x,x',x'',y\}$ with $x,x',x''\in V_i$ and $y\in V_{i+1}$ contains the three edges $\{x,x',y\}$, $\{x,x'',y\}$ and $\{x',x'',y\}$;
\end{itemize}
in both cases the fourth face is not an edge.

Since the chromatic index of $K_q$ is at most $q$ and $|V_2|=q$, there is a proper edge-coloring $\varphi$ of the complete graph $K_q[V_0]$ on $V_0$ with colors in $V_2$. Let
\[
D':=\bigl\{\{a,a',b,\varphi(aa')\}:\ a,a'\in V_0 \text{ and } b\in V_1\bigr\}.
\]
Each edge with two points in $V_0$ lies in a member of $D'$, and so does each edge $\{a,b,c\}$ that meets every part, unless no edge of $K_q[V_0]$ incident with $a$ has color $c$. Since $a$ is incident with $q-1$ edges and there are $q$ colors, each $a$ misses at most one color, so at most $|V_0||V_1|=q^2$ edges of $\mathcal{H}$ are not covered by $D'$.

For $i\in\{1,2\}$, let $\mathcal{T}_i$ be a family of triples of $V_i$ that covers every pair of $V_i$ and has $\bigl(\frac13+o(1)\bigr)\binom q2$ members; for instance, a partial Steiner triple system on $V_i$ that leaves $O(q)$ pairs uncovered (Theorem~\ref{r:5.11}), together with one triple through each of these uncovered pairs. Let
\[
D'':=\bigl\{T\cup\{y\}:\ T\in\mathcal{T}_i,\ y\in V_{i+1} \text{ and } i\in\{1,2\}\bigr\}.
\]
Every edge with two points in $V_1$ or in $V_2$ lies in a member of $D''$.

Let $D_4$ consist of the members of $D'$ and $D''$, together with one $4$-set containing each of the at most $q^2$ edges not covered by $D'$. Then every edge of $\mathcal{H}$ lies in a member of $D_4$, so $D_4$ also satisfies condition~(ii) of Proposition~\ref{r:4.1}, and
\[
|D_4|\le q\binom q2+2q\Bigl(\frac13+o(1)\Bigr)\binom q2+q^2=\Bigl(\frac56+o(1)\Bigr)q^3=\frac13\,e(\mathcal{H})+o(n^3).
\]
Therefore, by Proposition~\ref{r:4.1},
\[
\gamma\bigl(F_4(S_n)\bigr)\le\binom{n}{3}-e(\mathcal{H})+ |D_4|=\binom{n}{3}-\frac23\,e(\mathcal{H})+o(n^3)=\Bigl(\frac{17}{27}+o(1)\Bigr)\binom{n}{3}.
\]
If $n$ is not a multiple of $3$, we apply the construction to $3\lfloor n/3\rfloor$ of the points; the triples that meet the at most two remaining points are not edges of $\mathcal{H}$, and the same computation gives the bound.
\end{proof}

When $q\equiv1,3\pmod 6$, the families $\mathcal{T}_1$ and $\mathcal{T}_2$ can be taken to be Steiner triple systems (Theorem~\ref{r:5.11}). Moreover, since $q$ is odd, each $a\in V_0$ misses exactly one color $c_a$, and the $q^2$ edges $\{a,b,c_a\}$ not covered by $D'$ can be covered two at a time by sets $\{a,b,b',c_a\}$ with $b,b'\in V_1$, using $q(q+1)/2$ sets. 
This gives $|D_4|=q\binom q2+\frac{2q}{3}\binom q2+q\frac{q+1}{2}$, and therefore
\[
\gamma\bigl(F_4(S_{3q})\bigr)\le\binom{3q}{3}-e(\mathcal{H})+|D_4|=\binom{3q}{3}-\frac{q(q-1)(10q+3)}{6},
\]
which is $819$ for $n=21$ and $1809$ for $n=27$, against $966$ and $2160$ for the book construction. Together with~\eqref{eq:6},
\[
\bigl(0.628-o(1)\bigr)\binom{n}{3}\ \le\ \gamma\bigl(F_4(S_n)\bigr)\ \le\ \Bigl(\frac{17}{27}+o(1)\Bigr)\binom{n}{3},
\]
with $\frac{17}{27}\approx0.6296$. The upper constant is $1-\frac23\cdot\frac59$ and the lower one is $1-\frac23\pi_4$, so they coincide if Turán's conjecture $\pi_4=\frac59$ holds.

\subsection{Comparison with the bounds of Fabila-Monroy and Gómez-Galicia}\label{sec:4.5}

In this subsection we compare our bounds with those of Fabila-Monroy and Gómez-Galicia~\cite{b23}, the only previous work on the domination number of $F_k(S_n)$ for general $k$. Their lower bound is $\frac1k\binom{n}{k-1}$~\cite[Theorem~3.2]{b23}. Since $k(2k-3)-2(k-1)^2=k-2>0$, the constant $\frac{2k-3}{(k-1)^2}$ of Corollary~\ref{r:4.9} is more than twice $\frac1k$ for every $k\ge3$; for $k=3$ and $k=4$, the bounds in~\eqref{eq:6} give $\frac34$ and $0.628$, against $\frac13$ and $\frac14$.

Their upper bound is $\bigl(c'_k+o(1)\bigr)\binom{n}{k-1}$ with $c'_k=1-\frac1q+\frac{1}{k+q-1}$, where $q$ is the smallest prime dividing $k-1$ (denoted $p$ in~\cite{b23}); when $k\equiv1\pmod 4$ it improves to $c'_k=\frac12+\frac{1}{k+3}$. Ours is $c_k=1-\frac12\bigl(\frac{k-2}{k-1}\bigr)^{k-2}$ (Theorem~\ref{r:4.5}). Since $\bigl(\frac{k-2}{k-1}\bigr)^{k-2}=\bigl(1+\frac{1}{k-2}\bigr)^{-(k-2)}$ decreases to $e^{-1}$, the constant $c_k$ increases with $k$ from $c_3=\frac34$ towards $1-\frac{1}{2e}\approx0.8161$, so $\frac34\le c_k<0.8161$ for every $k\ge3$.

The comparison depends on $q$. If $q=2$, then $c'_k\le\frac12+\frac{1}{k+1}\le\frac34\le c_k$, with equality only for $k=3$, so their bound is at least as good. If $q=3$, then $k\in\{4,10,16,\dots\}$ and $c'_k=\frac23+\frac{1}{k+2}$; for $k=4$ ours is better, since $c_4=\frac79<\frac56=c'_4$ (and Theorem~\ref{r:4.12} lowers the constant to $\frac{17}{27}$), and for $k\ge10$ theirs is better, since $c'_k\le\frac34<c_k$. If $q\ge7$, then $c'_k>\frac67>0.8161>c_k$, so ours is better for every such $k$. Finally, if $q=5$, then $c'_k=\frac45+\frac{1}{k+4}$ decreases with $k$ while $c_k$ increases, so the two cross once; the values of $k$ with $q=5$ are those with $k-1\in\{5,25,35,55,65,\dots\}$, and a direct computation gives $c_k<c'_k$ for $k\in\{6,26,36,56\}$ and $c_k>c'_k$ for $k=66$, hence for every larger such $k$. In particular, $c_k$ is bounded away from $1$ uniformly in $k$, whereas $c'_k\to1$ along the values of $k$ for which $k-1$ is prime.

The construction of~\cite{b23} and the book construction are of the same kind: both choose the omitted $(k-1)$-sets so that few $k$-sets become complete (in~\cite{b23}, none, since the omitted sets are those whose sum is $1$ modulo $q$). By Theorem~\ref{r:4.8}, no such choice gives a constant below $1-\frac{k-2}{k-1}\pi_k$. For $k=3$ this is $\frac34$, which both constructions attain. 
For $k=4$ it lies between $0.628$, by~\eqref{eq:6}, and the constant $\frac{17}{27}\approx0.6296$ of Theorem~\ref{r:4.12}, whose construction has no complete $4$-set at all; the case $k\ge5$ is discussed in Section~\ref{sec:7}.

\section{Packing}\label{sec:5}

In this section we prove Theorem~\ref{r:1.1}, which identifies the packing number of $F_k(S_n)$ with the maximum number $\mathcal{D}(k-1,k,n+1)$ of blocks of a partial $S(k-1,k,n+1)$ on the $n+1$ vertices of the star, equivalently with the maximum size $A(n+1,4,k)$ of a binary code of length $n+1$, constant weight $k$ and minimum distance $4$. After that we prove the Steiner bound (which is used in Section~\ref{sec:6}), we present some small cases, and we state R\"odl's theorem, which shows that the bound is attained asymptotically. In Subsection~\ref{sec:5.2} we give two explicit constructions of packings by means of abelian groups.

\subsection{Packings and partial Steiner systems}\label{sec:5.1}

As in Section~\ref{sec:2}, we identify each vertex of $F_k(S_n)$ with a $k$-subset of $V=\{0,1,\dots,n\}$. Recall that a family of $k$-subsets of a set of $v$ points is a partial $S(k-1,k,v)$ if and only if any two distinct blocks share at most $k-2$ points. For the translation into codes, we denote by $A(v,d,k)$ the maximum size of a binary code of length $v$, constant weight $k$ and minimum Hamming distance $d$. Since the Hamming distance between the characteristic vectors of two blocks $B,B'$ is $|B\mathbin{\triangle} B'|=2(k-|B\cap B'|)$, the condition $|B\cap B'|\le k-2$ is equivalent to this distance being at least $4$, and so
\begin{equation}\label{eq:7}
\mathcal{D}(k-1,k,v)=A(v,4,k).
\end{equation}

\begin{teoremarep}{1.1}
Let $1\le k\le n$. A set $P$ of vertices of $F_k(S_n)$ is a packing if and only if $|A\cap B|\le k-2$ for all $A\neq B$ in $P$; that is, if and only if $P$ is a partial $S(k-1,k,n+1)$ on $V$. Consequently,
\[
\rho(F_k(S_n))=\mathcal{D}(k-1,k,n+1)=A(n+1,4,k).
\]
\end{teoremarep}
\begin{proof}
Let $A\neq B$ be vertices of $F_k(S_n)$, viewed as $k$-subsets of $V$, so that $|A\cap B|\le k-1$. If $A=\{0\}\cup X$ and $B=\{0\}\cup X'$ are centered, then $|A\cap B|=k-|X\setminus X'|$; if both are center-free, then $|A\cap B|=k-|A\setminus B|$; and if $A=\{0\}\cup X$ is centered and $B$ is center-free, then $|A\cap B|=k-1-|X\setminus B|$. Hence Lemma~\ref{r:3.6}, whose proof does not use the restriction $2\le k\le(n+1)/2$, reads
\[
d(A,B)=2\bigl(k-|A\cap B|\bigr)-\varepsilon,
\]
where $\varepsilon=1$ if exactly one of $A$ and $B$ is centered, and $\varepsilon=0$ otherwise. In particular $d(A,B)\le2$ if and only if $|A\cap B|=k-1$. Therefore $P$ is a packing if and only if $|A\cap B|\le k-2$ for all $A\neq B$ in $P$, which is the definition of a partial $S(k-1,k,n+1)$ on $V$. The last equality is~\eqref{eq:7}.
\end{proof}

The condition in Theorem~\ref{r:1.1} does not distinguish the center. Indeed, if $\pi$ is any permutation of $V$, then $\pi$ maps packings of $F_k(S_n)$ to packings of $F_k(S_n)$, even though in general $\pi$ does not induce an automorphism of $F_k(S_n)$. Therefore, when constructing packings we may label the $n+1$ points in any way, without worrying about which of them is the center; in Subsection~\ref{sec:5.2} the points will be the elements of an abelian group.

We now prove the Steiner bound, named after the designs that attain it with equality.

\begin{proposicion}\label{r:5.1}
For $2\le k\le n$,
\[
\rho(F_k(S_n))\le \frac{1}{k}\binom{n+1}{k-1},
\]
with equality if and only if a Steiner system $S(k-1,k,n+1)$ exists. In that case the blocks of the system form a maximum packing of $F_k(S_n)$.
\end{proposicion}

\begin{proof}
Let $P$ be a packing of $F_k(S_n)$. Each $A\in P$ contains $k$ subsets of size $k-1$, and by Theorem~\ref{r:1.1} none of them is contained in another element of $P$. Counting the pairs $(X,A)$ with $X\subset A\in P$ and $|X|=k-1$ we obtain $k|P|\le\binom{n+1}{k-1}$. Equality holds if and only if each $(k-1)$-subset of $V$ lies in exactly one element of $P$, that is, if $P$ is an $S(k-1,k,n+1)$.
\end{proof}

The following small cases illustrate Theorem~\ref{r:1.1} and Proposition~\ref{r:5.1}; we write $abc$ for $\{a,b,c\}$.

\begin{ejemplo}\label{r:5.3}\leavevmode
\begin{enumerate}[label=(\alph*)]
\item For $k=2$ the bound is $\frac{n+1}{2}$, and a packing is a family of pairwise disjoint pairs of $V$, so $\rho(F_2(S_n))=\lfloor\frac{n+1}{2}\rfloor$; the bound is attained exactly when $n$ is odd.
\item For $F_3(S_6)$ the bound is $\frac{1}{3}\binom{7}{2}=7$, attained by the Fano plane $\{013,124,235,346,045,156,026\}$, formed by the triples $\{i,i+1,i+3\}$ modulo $7$. Thus $\rho(F_3(S_6))=7$.
\item For $F_3(S_7)$ the bound is $\frac{1}{3}\binom{8}{2}=\frac{28}{3}$, which is not an integer, so no $S(2,3,8)$ exists. The Johnson bound (Proposition~\ref{r:5.2}) lowers it to $8$, and the eight triples $016,\allowbreak 025,\allowbreak 034,\allowbreak 124,\allowbreak 356,\allowbreak 157,\allowbreak 467,\allowbreak 237$ 
pairwise share at most one point, so $\rho(F_3(S_7))=8$.
\item For $F_6(S_{11})$ the bound is $\frac{1}{6}\binom{12}{5}=132$, and it is attained by the Witt design $S(5,6,12)$ (see~\cite{b15}).
\end{enumerate}
\end{ejemplo}

For fixed $k$ the Steiner bound is attained asymptotically; in the language of $F_k(S_n)$, this is a theorem of R\"odl.

\begin{teorema}[Rödl~\cite{b50}]\label{r:5.4}
For fixed $k\ge 2$,
\[
\rho(F_k(S_n))=\bigl(1-o(1)\bigr)\,\frac{1}{k}\binom{n+1}{k-1}\qquad(n\to\infty).
\]
\end{teorema}

By Proposition~\ref{r:5.1}, $\rho(F_k(S_n))=\frac{1}{k}\binom{n+1}{k-1}$ if and only if a Steiner system $S(k-1,k,n+1)$ exists. A necessary condition is that $(k-i)\mid\binom{n+1-i}{k-1-i}$ for $0\le i\le k-2$, and by Proposition~\ref{r:6.5}(i) these divisibility conditions amount to $\gcd(n-k+2,P_k)=1$, where $P_k$ is the product of the primes at most $k$. Keevash~\cite{b38} (see also~\cite{b28}) proved that, for fixed $k$, this condition is also sufficient when $n$ is large enough. For $k=3$ and $k=4$ it reads $n+1\equiv 1,3\pmod{6}$ and $n+1\equiv 2,4\pmod{6}$, respectively, and it is sufficient for every $n$ (Theorems~\ref{r:5.11} and~\ref{r:5.15}). For $k=5$ it is not: $17$ satisfies the divisibility conditions, but there is no $S(4,5,17)$~\cite{b46}.

Appendix~\ref{sec:B} contains the Johnson bound, the exact values of $\rho(F_k(S_n))$ for $k=3,4$ and the first open cases, all of which come from the literature through Theorem~\ref{r:1.1}. 

\subsection{Constructions with abelian groups}\label{sec:5.2}

We now give two constructions of packings that work for every $k$. As observed after Theorem~\ref{r:1.1}, it suffices to give a partial $S(k-1,k,n+1)$ on an arbitrary set of $n+1$ points. We take as points the elements of a finite abelian group $\Gamma$, written additively, possibly together with some new points, and for $A\subseteq\Gamma$ we write $\sigma(A)=\sum_{a\in A}a$. For $c\in\Gamma$ and $1\le k<|\Gamma|$, the \emph{sum class} of $c$ is
\[
\mathcal{P}_c(\Gamma,k)=\{A\subseteq\Gamma:\ |A|=k,\ \sigma(A)=c\}.
\]

\begin{lema}[see~\cite{b29}]\label{r:5.6}
Every sum class $\mathcal{P}_c(\Gamma,k)$ is a partial $S(k-1,k,|\Gamma|)$.
\end{lema}

\begin{proof}
If $A,A'\in\mathcal{P}_c(\Gamma,k)$ contain the same $(k-1)$-subset $X$, then $A=X\cup\{a\}$, $A'=X\cup\{a'\}$ and $\sigma(X)+a=c=\sigma(X)+a'$; hence $a=a'$ and $A=A'$.
\end{proof}

Since the $|\Gamma|$ sum classes partition the $k$-subsets of $\Gamma$, one of them has at least $\binom{|\Gamma|}{k}/|\Gamma|$ elements. Graham and Sloane~\cite{b29} used this with $\Gamma=\mathbb{Z}_{n+1}$, and with Theorem~\ref{r:1.1} it gives the following bound.

\begin{corolario}[\cite{b29}]\label{r:5.7}
For $1\le k\le n$, $\rho(F_k(S_n))\ge\dfrac{1}{n+1}\dbinom{n+1}{k}$.
\end{corolario}

The size of a class depends on $\Gamma$ and on $c$, and Kløve~\cite{b40} determined its maximum over all abelian groups $\Gamma$ of order $v$ and all $c\in\Gamma$. Write $v=\prod_p p^{E_p}$, and let $\Gamma_v^*=\prod_p\mathbb{Z}_p^{E_p}$ be the abelian group of order $v$ whose Sylow subgroups are elementary abelian.

\begin{teorema}[Kløve~\cite{b40}]\label{r:5.8}
Let $1\le k<v$, $g=\gcd(v,k)$, and let $\mathcal{R}$ be the set of squarefree divisors of $g$. For $r\in\mathcal{R}$ let $w_k(r)=\prod_p(p^{E_p}-1)$, where the product runs over the primes $p\mid r$, except $p=2$ when $k\equiv 2\pmod{4}$. Let
\[
K(v,k)=\frac{1}{v}\sum_{r\in\mathcal{R}}\binom{v/r}{k/r}w_k(r).
\]
Then $|\mathcal{P}_c(\Gamma,k)|\le K(v,k)$ for every abelian group $\Gamma$ of order $v$ and every $c\in\Gamma$, and equality holds for $\Gamma=\Gamma_v^*$ with $c=0$ if $k\not\equiv 2\pmod{4}$ or $v$ is odd, and with $c$ an element of order $2$ of the factor $\mathbb{Z}_2^{E_2}$ if $k\equiv 2\pmod{4}$ and $v$ is even. Consequently, $\rho(F_k(S_n))\ge K(n+1,k)$.
\end{teorema}

A sum class need not be a maximum packing, and it can be enlarged with new points. If $\mathcal{B}$ is a partial $S(k-1,k,v)$ on a set $X$, its \emph{leave} $L(\mathcal{B})$ is the family of the $(k-1)$-subsets of $X$ that lie in no block of $\mathcal{B}$.

\begin{lema}\label{r:5.9}
Let $\Gamma$ be an abelian group of order $m$, $c\in\Gamma$, $\mathcal{P}=\mathcal{P}_c(\Gamma,k)$, and let $I$ be a set of $j$ new points. Let $\mathcal{Q}$ be a family of $k$-subsets of $\Gamma\cup I$ such that
\begin{enumerate}[label=(\alph*)]
\item $\mathcal{Q}$ is a partial $S(k-1,k,m+j)$, and
\item every $(k-1)$-subset of $\Gamma$ contained in a block of $\mathcal{Q}$ belongs to $L(\mathcal{P})$.
\end{enumerate}
Then $\mathcal{P}\cup\mathcal{Q}$ is a partial $S(k-1,k,m+j)$, and therefore $\rho(F_k(S_{m+j-1}))\ge|\mathcal{P}|+|\mathcal{Q}|$.
\end{lema}
\begin{proof}
Let $X$ be a $(k-1)$-subset of $\Gamma\cup I$. By Lemma~\ref{r:5.6} and~(a), $X$ lies in at most one block of $\mathcal{P}$ and in at most one block of $\mathcal{Q}$, so it suffices to show that it does not lie in a block of each. If $X$ lies in a block of $\mathcal{P}$, then $X\subseteq\Gamma$ and $X\notin L(\mathcal{P})$, so by~(b) it lies in no block of $\mathcal{Q}$. In particular, no block belongs to both $\mathcal{P}$ and $\mathcal{Q}$, so $|\mathcal{P}\cup\mathcal{Q}|=|\mathcal{P}|+|\mathcal{Q}|$, and the bound follows from Theorem~\ref{r:1.1}.
\end{proof}

We call $\mathcal{P}_c(\Gamma,k)\cup\mathcal{Q}$ as in Lemma~\ref{r:5.9} a \emph{$\Gamma+j$ points construction}. A block of $\mathcal{Q}$ either lies in $\Gamma$, and then all its $(k-1)$-subsets belong to $L(\mathcal{P})$, or contains $r\ge1$ new points; if $r=1$ it is $X\cup\{\infty\}$ with $X\in L(\mathcal{P})$,\footnote{As usual in design theory, we denote the new points by $\infty$, or by $\infty_1,\dots,\infty_j$ when there are several.} and if $r\ge2$ condition (b) holds automatically. Neither construction is optimal in general, and each is better than the other for some parameters; in Appendix~\ref{sec:C} we compare them on small examples and, for $k=3$, with the exact value.

\section{Domination and packing}\label{sec:6}

In every graph we have $\rho\le\gamma$ (Section~\ref{sec:2}). In this section we prove Theorem~\ref{r:1.3}, which determines exactly when equality holds for $F_k(S_n)$, namely only in the self-complementary range $n=2k-1$, and there if and only if a Steiner system $S(k-2,k-1,2k-1)$ exists, which by a theorem of Hammond and Smith~\cite{b31} is the same as a perfect 1-code in the odd graph $O_k$. The ingredients are a common bound for the two parameters and the characterization of the perfect dominating sets of $F_k(S_{2k-1})$ by means of Steiner systems (Subsection~\ref{sec:6.2}). 

Let $2\le k\le(n+1)/2$. By Corollary~\ref{r:3.2}, the maximum degree of $F_k(S_n)$ is $n-k+1$, so each vertex dominates at most $n-k+2$ vertices, and every dominating set has at least $|V(F_k(S_n))|/(n-k+2)$ elements. Together with the Steiner bound (Proposition~\ref{r:5.1}), this gives
\begin{equation}\label{eq:common}
\rho\bigl(F_k(S_n)\bigr)\;\le\;\frac{1}{k}\binom{n+1}{k-1}\;=\;\frac{|V(F_k(S_n))|}{n-k+2}\;\le\;\gamma\bigl(F_k(S_n)\bigr),
\end{equation}
where the equality holds because $|V(F_k(S_n))|=\binom{n+1}{k}=\frac{n-k+2}{k}\binom{n+1}{k-1}$. 
Thus the two classical bounds, the Steiner bound for packing and the volume bound for domination, are the same number.
In general neither is attained; for fixed $k\ge3$, Corollary~\ref{r:4.9} gives $\gamma(F_k(S_n))\ge\bigl(\frac{2k-3}{(k-1)^2}-o(1)\bigr)\binom{n}{k-1}$, so that $\gamma(F_k(S_n))\ge(2-o(1))\,\rho(F_k(S_n))$ as $n\to\infty$.

\subsection{Perfect domination and Steiner systems}\label{sec:6.2}

In this subsection we show that, in the self-complementary range $n=2k-1$, the existence of a perfect dominating set is equivalent to the existence of a Steiner system $S(k-2,k-1,2k-1)$.

For $n=2k-1$ the graph $F_k(S_{2k-1})$ is $k$-regular (Corollary~\ref{r:3.2}) and has $\binom{2k}{k}$ vertices. 
By~\eqref{eq:common} with $n=2k-1$, since $\frac1k\binom{2k}{k-1}=\frac{1}{k+1}\binom{2k}{k}$,
\begin{equation}\label{eq:9}
\rho\bigl(F_k(S_{2k-1})\bigr)\;\le\;C_k\;\le\;\gamma\bigl(F_k(S_{2k-1})\bigr),\qquad C_k=\frac{1}{k+1}\binom{2k}{k},
\end{equation}
where $C_k$ is the $k$-th Catalan number. Since the graph is regular, a perfect dominating set (a \emph{perfect 1-code}, in the terminology of coding theory~\cite{b8}) if it exists, has exactly $C_k$ elements, and $\gamma(F_k(S_{2k-1}))=C_k$ if and only if $F_k(S_{2k-1})$ has a perfect dominating set.

For $0\le i\le t$, the number of blocks of an $S(t,r,v)$ containing a fixed $i$-subset $I$ is
\begin{equation}\label{eq:10}
\lambda_i=\binom{v-i}{t-i}\Big/\binom{r-i}{t-i},
\end{equation}
independent of $I$, by counting the pairs $(T,B)$ with $I\subseteq T\subseteq B$, $|T|=t$ and $B$ a block (see \cite[II.5]{b15}); in particular the system has $\lambda_0$ blocks, and the integrality of all the $\lambda_i$ is necessary for it to exist.

We write $[n]=[2k-1]$ for the set of leaves. The proofs of this subsection use the following property of the systems $S(k-2,k-1,2k-1)$, whose proof, an inclusion--exclusion count over the numbers $\lambda_i$, is given in Appendix~\ref{sec:C}.

\begin{lema}\label{r:6.2}
Let $k\ge3$ and let $\mathcal{S}$ be a Steiner system $S(k-2,k-1,n)$ on $[n]=[2k-1]$. Then:
\begin{enumerate}[label=(\roman*)]
\item $k$ is even;
\item every $k$-subset of $[n]$ contains at most one block;
\item for every $(k-1)$-subset $X$ of $[n]$, the number of blocks disjoint from $X$ is $0$ if $X\in\mathcal{S}$ and $1$ if $X\notin\mathcal{S}$;
\item a $k$-subset $A$ of $[n]$ contains no block if and only if $[n]\setminus A$ is a block.
\end{enumerate}
\end{lema}

In the odd graph $O_k$, whose vertices are the $(k-1)$-subsets of $[n]$ and whose edges join disjoint subsets, Lemma~\ref{r:6.2}(iii) says that every vertex is either a block or adjacent to exactly one block, that is, the blocks form a perfect 1-code of $O_k$.

The next theorem combines, in the language of token graphs, two results of Hammond and Smith~\cite{b31} on perfect codes in the odd graphs, in the framework introduced by Biggs~\cite{b8}: their Theorem~3, by which $O_k$ has a perfect 1-code if and only if a Steiner system $S(k-2,k-1,2k-1)$ exists, and their Theorem~2, by which $O_k$ has a perfect $e$-code if and only if $2O_k$ does. By~\eqref{eq:2}, $F_k(S_{2k-1})$ is the bipartite double $2O_k$, a $k$-set $A$ corresponding to the $(k-1)$-set $[n]\setminus A$, so that $X\subset A$ if and only if $X$ and $[n]\setminus A$ are disjoint; the set $D$ below is the lift to $2O_k$ of the perfect 1-code formed by the blocks of $\mathcal{S}$. 
Our proof is direct and short, avoids the argument on which the converse of their Theorem~2 rests (their Theorem~1, whose proof uses Biggs' generalization of Lloyd's theorem~\cite{b8}), and gives the perfect dominating set explicitly in terms of the Steiner system. 

\begin{teorema}[{see~\cite[Theorems~2 and~3]{b31}}]\label{r:6.3}
Let $k\ge2$. The graph $F_k(S_{2k-1})$ has a perfect dominating set if and only if there exists a Steiner system $S(k-2,k-1,2k-1)$. Moreover, if $\mathcal{S}$ is such a system, then
\[
D=\bigl\{\,\{0\}\cup X:X\in\mathcal{S}\,\bigr\}\;\cup\;\bigl\{\,[2k-1]\setminus X:X\in\mathcal{S}\,\bigr\}
\]
is a perfect dominating set of size $C_k$.
\end{teorema}
\begin{proof}
We work in $B(2k-1;k-1,k)$ by means of Theorem~\ref{r:3.1}, with $[n]=[2k-1]$ the set of leaves, and for a condition $P$ we write $[\,P\,]=1$ if $P$ holds and $[\,P\,]=0$ otherwise. For $k=2$ both conditions hold: the system $S(0,1,3)$ consists of a single block $\{a\}$, and $D=\{\{0,a\},[3]\setminus\{a\}\}$ is a perfect dominating set of $F_2(S_3)=C_6$ with $C_2=2$ elements. In what follows we assume that $k\ge3$.

Suppose first that $\mathcal{S}$ is an $S(k-2,k-1,2k-1)$, and define $D$ as in the statement. Since $|[n]\setminus X|=k$ for every block $X$, the second part of $D$ consists of center-free vertices, so $|D|=2|\mathcal{S}|$. We show that every vertex is dominated exactly once.

Let $A$ be a center-free vertex, that is, a $k$-subset of $[n]$. Its centered neighbors in $D$ are the sets $\{0\}\cup X$ with $X\in\mathcal{S}$ and $X\subset A$, and by Lemma~\ref{r:6.2}(ii) there is at most one of them; $A$ has no center-free neighbors, since center-free vertices are pairwise nonadjacent. If $A$ contains a block, then $[n]\setminus A$ is not a block by Lemma~\ref{r:6.2}(iv), so $A\notin D$, and $A$ is dominated exactly once, by a centered vertex. If $A$ contains no block, then $[n]\setminus A$ is a block by Lemma~\ref{r:6.2}(iv), so $A\in D$ and $A$ has no neighbor in $D$; again it is dominated exactly once.

Now let $\{0\}\cup X$ be a centered vertex, with $X$ a $(k-1)$-subset of $[n]$. It lies in $D$ if and only if $X\in\mathcal{S}$, and its neighbors in $D$ are the sets $[n]\setminus Y$ with $Y\in\mathcal{S}$ and $X\subset[n]\setminus Y$, that is, with $Y$ a block disjoint from $X$. By Lemma~\ref{r:6.2}(iii), if $X\in\mathcal{S}$ then $\{0\}\cup X$ lies in $D$ and has no neighbor in $D$, and if $X\notin\mathcal{S}$ then it does not lie in $D$ and has exactly one neighbor in $D$. In both cases it is dominated exactly once.

Therefore $D$ is a perfect dominating set, and since $F_k(S_{2k-1})$ is $k$-regular with $\binom{2k}{k}$ vertices, $|D|=\binom{2k}{k}/(k+1)=C_k$.

Conversely, let $D$ be a perfect dominating set of $F_k(S_{2k-1})$, and let
\[
D_{k-1}=\{\,X:\{0\}\cup X\in D\,\},\qquad D_k=\{\,A\in D:A\subseteq[n]\,\}
\]
be its centered and center-free parts, viewed as families of $(k-1)$- and of $k$-subsets of $[n]$. Since every vertex is dominated exactly once,
\begin{align}
[\,A\in D_k\,]+\bigl|\{X\in D_{k-1}:X\subset A\}\bigr|&=1\quad\text{for every $k$-set }A,\label{eq:11}\\
[\,X\in D_{k-1}\,]+\bigl|\{A\in D_k:X\subset A\}\bigr|&=1\quad\text{for every $(k-1)$-set }X.\label{eq:12}
\end{align}
Each $(k-1)$-set lies in $k$ sets of size $k$, and each $k$-set contains $k$ sets of size $k-1$, so summing \eqref{eq:11} over all $k$-sets and \eqref{eq:12} over all $(k-1)$-sets gives
\[
|D_k|+k|D_{k-1}|=\binom{2k-1}{k},\qquad |D_{k-1}|+k|D_k|=\binom{2k-1}{k-1}.
\]
The right-hand sides are equal, so $(k-1)\bigl(|D_{k-1}|-|D_k|\bigr)=0$; since $|D|=C_k$ by regularity, $|D_{k-1}|=|D_k|=C_k/2$.

We claim that $D_{k-1}$ is an $S(k-2,k-1,2k-1)$. Let $W$ be a $(k-2)$-subset of $[n]$. If $W$ were contained in two distinct elements $W\cup\{x_1\}$ and $W\cup\{x_2\}$ of $D_{k-1}$, the $k$-set $W\cup\{x_1,x_2\}$ would contain both, and the left-hand side of \eqref{eq:11} would be at least $2$; so $W$ lies in at most one element of $D_{k-1}$. On the other hand, each element of $D_{k-1}$ contains $k-1$ sets of size $k-2$, so, using $\binom{2k}{k}=2\binom{2k-1}{k-1}$, the number of pairs $(W,X)$ with $W\subset X\in D_{k-1}$ is
\[
(k-1)\,|D_{k-1}|=\frac{k-1}{2}\,C_k=\frac{k-1}{k+1}\binom{2k-1}{k-1}=\binom{2k-1}{k-2},
\]
the number of $(k-2)$-subsets of $[n]$. Hence every $(k-2)$-set lies in exactly one element of $D_{k-1}$.
\end{proof}

Theorem~\ref{r:6.3} turns the question into the existence of a Steiner system, and the divisibility conditions~\eqref{eq:10} give a first necessary condition. For an $S(k-1,k,v)$, that is, with $t=k-1$ and $r=k$, we have $\binom{k-i}{k-1-i}=k-i$, so these conditions read
\begin{equation}\label{eq:adm}
(k-i)\;\Big|\;\binom{v-i}{k-1-i}\qquad(0\le i\le k-2),
\end{equation}
and we call an integer $v\ge k$ \emph{$k$-admissible} if they hold.

The $k-1$ conditions~\eqref{eq:adm} reduce to a single one. Let $P_k=\prod_{p\le k}p$ be the product of the primes less than or equal to $k$. The necessity of ``$k+1$ prime'' for the existence of an $S(k-2,k-1,2k-1)$, part~(ii) below, is Lemma~3 of~\cite{b31}; part~(i) places it in a criterion valid for every $v$.

\begin{proposicion}\label{r:6.5}
Let $k\ge2$.
\begin{enumerate}[label=(\roman*)]
\item An integer $v\ge k$ is $k$-admissible if and only if $\gcd(v-k+1,\,P_k)=1$.
\item The Steiner system $S(k-2,k-1,2k-1)$ satisfies all the necessary divisibility conditions if and only if $2k$ is $k$-admissible, if and only if $k+1$ is prime.
\end{enumerate}
\end{proposicion}

The proof is an elementary divisibility argument, and we defer it to Appendix~\ref{sec:C}.

\begin{corolario}\label{r:6.6}
For $k\in\{2,4,6\}$ the graph $F_k(S_{2k-1})$ has a perfect dominating set, obtained from the trivial system $S(0,1,3)$, the Fano plane $S(2,3,7)$ and the Witt design $S(4,5,11)$, respectively, and $\gamma(F_k(S_{2k-1}))=C_k$. For $k\in\{3,5,7\}$ it has none, and $\gamma(F_k(S_{2k-1}))>C_k$.
\end{corolario}

\begin{proof}
The three systems exist, and for $k\in\{3,5,7\}$ no $S(k-2,k-1,2k-1)$ exists by Proposition~\ref{r:6.5}(ii), since $k+1$ is not prime. The rest follows from Theorem~\ref{r:6.3} and the fact that $\gamma(F_k(S_{2k-1}))=C_k$ if and only if $F_k(S_{2k-1})$ has a perfect dominating set.
\end{proof}

The sizes are $C_2=2$, $C_4=14$ and $C_6=132$; in particular, $\gamma(F_6(S_{11}))=132$ follows from the Witt design without any computation. By Proposition~\ref{r:6.5}(ii), $F_k(S_{2k-1})$ can have a perfect dominating set only if $k+1$ is prime, that is, for $k\in\{2,4,6,10,12,16,\dots\}$, and the first open case is $k=10$, which corresponds to the system $S(8,9,19)$. To the best of our knowledge, no Steiner system $S(t,r,v)$ with $t\ge6$ is explicitly known (see~\cite[II.5]{b15}). Keevash's theorem~\cite{b38} guarantees that they exist when $v$ is large in terms of $t$ and $r$, but it does not apply to the systems $S(k-2,k-1,2k-1)$, where $v=2t+3$ (Problem~\ref{r:7.5}).

\subsection{When packing and domination coincide}\label{sec:6.3}

We now prove Theorem~\ref{r:1.3}, which we restate for easy reference.

\begin{teoremarep}{1.3}
Let $2\le k\le n-1$. Then $\rho(F_k(S_n))=\gamma(F_k(S_n))$ if and only if $n=2k-1$ and there exists a Steiner system $S(k-2,k-1,2k-1)$, in which case $\rho(F_k(S_n))=\gamma(F_k(S_n))=C_k=\frac{1}{k+1}\binom{2k}{k}$, the $k$-th Catalan number.
\end{teoremarep}

\begin{proof}
By~\eqref{eq:1}, and since $n=2k-1$ forces $n+1-k=k$, we may assume $2\le k\le(n+1)/2$. Write $V=V(F_k(S_n))$, $\rho=\rho(F_k(S_n))$, $\gamma=\gamma(F_k(S_n))$ and $\Delta=\Delta(F_k(S_n))$; by Corollary~\ref{r:3.2}, $\Delta=n-k+1$.

Suppose that $\rho=\gamma$, and let $D$ be a dominating set with $|D|=\gamma$. By~\eqref{eq:common}, $|D|=|V|/(\Delta+1)$, so
\[
|V|\;\le\;\sum_{A\in D}|N[A]|\;\le\;|D|\,(\Delta+1)\;=\;|V|.
\]
Hence equality holds throughout, so every vertex of $D$ has degree $\Delta$ and the closed neighborhoods of the vertices of $D$ are pairwise disjoint; that is, $D$ is a perfect dominating set. If $n>2k-1$, then $\Delta>k$, so by Corollary~\ref{r:3.2} every vertex of $D$ is centered. Since centered vertices are pairwise nonadjacent, a centered vertex outside $D$ would not be dominated, so $D$ consists of all centered vertices. But then each center-free vertex, whose $k\ge2$ neighbors are all centered, lies in at least two closed neighborhoods of vertices of $D$, a contradiction. Therefore $n=2k-1$, and Theorem~\ref{r:6.3} gives a Steiner system $S(k-2,k-1,2k-1)$.

Conversely, if $n=2k-1$ and such a system exists, Theorem~\ref{r:6.3} gives a perfect dominating set $D$ with $|D|=C_k$. Since $D$ is both a packing and a dominating set, $\gamma\le|D|\le\rho\le\gamma$, so $\rho=\gamma=C_k$.
\end{proof}

A perfect dominating set is at the same time a packing and a dominating set, so its existence forces $\rho(F_k(S_n))=\gamma(F_k(S_n))$. We thus obtain the following.

\begin{corolario}\label{r:6.8}
Let $2\le k\le n-1$. The graph $F_k(S_n)$ has a perfect dominating set if and only if $n=2k-1$ and a Steiner system $S(k-2,k-1,2k-1)$ exists.
\end{corolario}

\begin{proof}
If $D$ is a perfect dominating set of $F_k(S_n)$, then $\gamma(F_k(S_n))\le|D|\le\rho(F_k(S_n))\le\gamma(F_k(S_n))$, so the two parameters are equal, and Theorem~\ref{r:1.3} gives $n=2k-1$ and the existence of the system. The converse is Theorem~\ref{r:6.3}.
\end{proof}

For $n=2k-1$ the Steiner bound is $\frac1k\binom{2k}{k-1}=C_k$, so by Proposition~\ref{r:5.1} the packing number of $F_k(S_{2k-1})$ equals $C_k$ exactly when an $S(k-1,k,2k)$ exists; moreover, a packing with $C_k$ elements is a perfect dominating set, since its $C_k$ closed neighborhoods are disjoint, have $k+1$ vertices each, and $(k+1)C_k=\binom{2k}{k}$. Together with Theorem~\ref{r:1.3}, this shows that an $S(k-1,k,2k)$ exists if and only if an $S(k-2,k-1,2k-1)$ does, in agreement with design theory: the derived design of an $S(k-1,k,2k)$ at a point is an $S(k-2,k-1,2k-1)$, and conversely, since $k-2$ is even by Lemma~\ref{r:6.2}(i), an $S(k-2,k-1,2k-1)$ extends by complementation to an $S(k-1,k,2k)$, by a theorem of Alltop~\cite{b4}.

\section{Concluding remarks and open problems}\label{sec:7}

All the results of Section~\ref{sec:3} and Appendix~\ref{sec:A} come from the identification $F_k(S_n)\cong B(n;k-1,k)$. Two parameters that it does not settle are the crossing number and the genus. For $k=2$ they are those of $K_n$ (Proposition~\ref{r:3.5}), because $F_2(S_n)$ is a subdivision of $K_n$; as noted in Appendix~\ref{sec:A1}, $F_k(S_n)$ is a subdivision only when $k=2$, and for $k\ge3$ we only know that it is not planar (Proposition~\ref{r:3.8}(v)).

\begin{problema}\label{r:7.1}
Determine the crossing number and the genus of $F_k(S_n)$ for $k\ge3$.
\end{problema}

The remaining problems concern domination and packing. For the domination number, Theorem~\ref{r:1.2} gives the exact constant $\frac34$ for $k=3$; for $k=4$ it gives, with~\eqref{eq:6}, the constants $0.628$ and $\frac{17}{27}\approx0.6296$, and for $k\ge5$ it leaves a gap between $1-\frac{k-2}{k-1}\pi_k$ and $1-\frac12\bigl(\frac{k-2}{k-1}\bigr)^{k-2}$. For $k\ge5$ the upper bound comes from one family of hypergraphs, the books of Construction~\ref{r:4.3}, and the lower bound from the removal lemma. For a $K_k^{(k-1)}$-free hypergraph $\mathcal{H}$ we have $d_k(\mathcal{H})\ge e(\mathcal{H})/(k-1)$ by Lemma~\ref{r:4.2}. If, in addition, all but $o(n^{k-1})$ edges of $\mathcal{H}$ can be covered by $k$-sets that contain exactly $k-1$ edges each and pairwise share no edge, then $d_k(\mathcal{H})=\bigl(1+o(1)\bigr)e(\mathcal{H})/(k-1)$, and Proposition~\ref{r:4.1} gives $\gamma(F_k(S_n))\le\binom{n}{k-1}-\frac{k-2}{k-1}\,e(\mathcal{H})+o(n^{k-1})$. For $k=3$ the balanced complete bipartite graph has this property, since its edges split into paths of length two, and it gives the constant $\frac34$ again; for $k=4$, the proof of Theorem~\ref{r:4.12} shows that Turán's hypergraph has it, and it gives $\frac{17}{27}$. For $k=3$, moreover, integer programming shows that the books are optimal for all $3\le n\le10$ (Subsection~\ref{sec:4.4}), and we conjecture that this persists.

\begin{conjetura}\label{r:7.2}
For every $n\ge3$, $\gamma(F_3(S_n))=\binom n2-\bigl\lfloor(n+1)^2/8\bigr\rfloor$.
\end{conjetura}

By Corollary~\ref{r:4.11}, Conjecture~\ref{r:7.2} is asymptotically true, and what remains open is the exact value. For $k\ge4$, the discussion above suggests that the lower bound of Theorem~\ref{r:4.8} is the right one.

\begin{conjetura}\label{r:7.3}
For every $k\ge3$, $\displaystyle\lim_{n\to\infty}\gamma\bigl(F_k(S_n)\bigr)\Big/\binom{n}{k-1}=1-\frac{k-2}{k-1}\,\pi_k$.
\end{conjetura}

For $k=3$, Conjecture~\ref{r:7.3} is Corollary~\ref{r:4.11}, and for $k=4$ it follows from Turán's conjecture $\pi_4=\frac59$, by Theorem~\ref{r:4.12}.

The packing and the domination numbers are compared in Section~\ref{sec:6} only at equality. For fixed $k\ge3$, \eqref{eq:common} and Corollary~\ref{r:4.9} give $\gamma(F_k(S_n))\ge(2-o(1))\,\rho(F_k(S_n))$ as $n\to\infty$, but when equality fails for a fixed $n$ we do not know how large the difference is, not even in the self-complementary range. For instance, $\rho(F_3(S_5))=4<C_3=5<6=\gamma(F_3(S_5))$, and outside that range $\rho$ can attain the common bound~\eqref{eq:common} without equality: $\rho(F_3(S_8))=12=\frac13\binom92$, since an $S(2,3,9)$ exists, whereas $\gamma(F_3(S_8))=18$.

\begin{problema}\label{r:7.4}
Determine $\gamma(F_k(S_{2k-1}))-C_k$ when no Steiner system $S(k-2,k-1,2k-1)$ exists, in particular for odd $k$, where no such system exists by Lemma~\ref{r:6.2}(i).
\end{problema}

Finally, Theorem~\ref{r:6.3} reduces the existence of a perfect dominating set of $F_k(S_{2k-1})$ to that of an $S(k-2,k-1,2k-1)$, which requires $k+1$ to be prime (Proposition~\ref{r:6.5}) and is known only for $k\in\{2,4,6\}$. Hammond and Smith~\cite{b31} observed that if, as was then conjectured, no $t$-designs with $t>5$ existed, then $O_4$ and $O_6$ would be, apart from trivial cases, the only odd graphs with a perfect 1-code; in our language, $F_2(S_3)=C_6$, $F_4(S_7)$ and $F_6(S_{11})$ would be the only token graphs of stars with $2\le k\le n-1$ and $\rho=\gamma$. That conjecture is false, and for Steiner systems, which are the relevant case here, Keevash~\cite{b38} proved that $S(t,r,v)$ exists for every $t$ whenever $v$ is large enough and the divisibility conditions hold. His theorem does not apply to the family $S(k-2,k-1,2k-1)$, where $v=2t+3$, so whether these three are the only ones remains open.

\begin{problema}\label{r:7.5}
Does $F_{10}(S_{19})$ have a perfect dominating set? Equivalently, does a Steiner system $S(8,9,19)$ exist?
\end{problema}

\section*{Declaration on the use of AI tools}
During the preparation of this work, the authors used Claude (Anthropic) to discuss their approach and techniques as they developed them into the detailed strategies and constructions presented here. It was also used to verify computationally the values of the parameters of $F_k(S_n)$ for small values of $n$ and $k$, to search for references, and to edit the text. The authors posed the problems, checked every proof and computation, verified each reference against its source, and take full responsibility for the content of the paper.

\bibliographystyle{abbrv}
\bibliography{bibliography-star}

\appendix

\section{Further parameters of $F_k(S_n)$}\label{sec:A}

In this appendix we collect the parameters of $F_k(S_n)$ that are either known, as special cases of results on the token graphs of wider classes of graphs or on the doubled Johnson graphs, or almost immediate from Theorem~\ref{r:3.1}. For the parameters whose proofs lie outside the scope of this paper, the genus and the crossing number of $K_n$, we state the results with precise references; for all the others we give proofs in the language of the inclusion graph $B(n;k-1,k)$. Throughout we assume $2\le k\le(n+1)/2$, which by~\eqref{eq:1} is no loss of generality, except in Subsection~\ref{sec:A1}, where $k=2$.

\subsection{The case $k=2$}\label{sec:A1}

By Corollary~\ref{r:3.4}, $F_2(S_n)$ is the total subdivision of $K_n$; the graph on the left of Figure~\ref{fig:1} is the case $n=5$. This is special to $k=2$, since by Corollary~\ref{r:3.2} the center-free side of $F_k(S_n)$ has degree $k$, so that $F_k(S_n)$ is a subdivision only for $k=2$. Since subdividing edges changes neither planarity, nor the crossing number, nor the genus, these three parameters transfer from $K_n$ to $F_2(S_n)$ at once.

\begin{proposicion}\label{r:3.5}
The graph $F_2(S_n)$ is planar if and only if $n\le4$, and
\[
\cro\bigl(F_2(S_n)\bigr)=\cro(K_n),\qquad \gen\bigl(F_2(S_n)\bigr)=\gen(K_n)=\left\lceil\frac{(n-3)(n-4)}{12}\right\rceil.
\]
\end{proposicion}

\begin{proof}
Subdividing edges preserves planarity, the crossing number and the genus, and $K_n$ is planar exactly for $n\le4$. The value of $\gen(K_n)$ was determined by Ringel and Youngs~\cite{b49}, thereby settling the map-color conjecture of Heawood~\cite{b33}.
\end{proof}

For $k=2$, where $F_1(S_n)=S_n$, Theorem~\ref{r:3.3} says that $F_2(S_{n+1})$ arises from $F_2(S_n)$ by adding the centered configuration $\{0,n+1\}$ and the $n$ center-free configurations $\{i,n+1\}$, each joined to $\{0,i\}$ and to $\{0,n+1\}$; through Corollary~\ref{r:3.4} this is the step $K_n\to K_{n+1}$ with each of the $n$ new edges subdivided.

The equality $\cro(F_2(S_n))=\cro(K_n)$ transfers to $F_2(S_n)$ what is known about $\cro(K_n)$. The Harary--Hill conjecture states that
\[
\cro(K_n)=Z(n)=\frac14\Bigl\lfloor\frac n2\Bigr\rfloor\Bigl\lfloor\frac{n-1}2\Bigr\rfloor\Bigl\lfloor\frac{n-2}2\Bigr\rfloor\Bigl\lfloor\frac{n-3}2\Bigr\rfloor
\]
for every $n\ge3$. It was known for $n\le14$~\cite{b47,b2}, and Balogh, Lidick\'y and Salazar~\cite{b6} proved that $\cro(K_n)\ge0.985\,Z(n)$ for large $n$. A proof for every $n$ was recently given in a preprint produced by an OpenAI model~\cite{b44}, with a Lean formalization, and Salazar~\cite{bSal} has written a detailed combinatorial exposition of it, together with an extension to rotation systems. By Proposition~\ref{r:3.5}, $\cro(F_2(S_n))=Z(n)$ for every $n\ge3$.

\subsection{Colorings, connectivity, independence and planarity}\label{sec:A2}

Parts \emph{(iii)}, \emph{(iv)} and \emph{(v)} of the next proposition are special cases of known results. Indeed, the equality $\kappa(F_k(T))=\lambda(F_k(T))=\delta(F_k(T))$ holds for the token graph of every tree $T$~\cite{b26}, and \emph{(iii)} is the case $T=S_n$; the independence number of $F_k(S_n)$ was determined by de Alba, Carballosa, Lea\~nos and Rivera~\cite{b19}, and the matching number in \emph{(iv)} then follows from K\H{o}nig's theorem; and the planar token graphs were characterized by Carballosa, Fabila-Monroy, Lea\~nos and Rivera~\cite{b14}, so \emph{(v)} is the case $G=S_n$ of their result. We prove them in the language of $B(n;k-1,k)$.
\begin{proposicion}\label{r:3.8}
For $2\le k\le(n+1)/2$ the graph $F:=F_k(S_n)$ satisfies:
\begin{enumerate}[leftmargin=*]
\item[(\textit{i})] $\chi(F)=2$ and $\chi'(F)=\Delta(F)=n-k+1$;
\item[(\textit{ii})] $\operatorname{girth}(F)=6$;
\item[(\textit{iii})] $\kappa(F)=\lambda(F)=\delta(F)=k$;
\item[(\textit{iv})] $\alpha(F)=\binom nk$ and $\nu(F)=\tau(F)=\binom{n}{k-1}$;
\item[(\textit{v})] $F$ is planar if and only if $n\le4$, which in this range forces $k=2$.
\end{enumerate}
\end{proposicion}

\begin{proof}
\emph{(i)} The graph is bipartite with at least one edge, so $\chi(F)=2$, and K\H{o}nig's edge-coloring theorem~\cite{bKonig16} gives $\chi'=\Delta$ in every bipartite graph, where $\Delta(F)=n-k+1$ by Corollary~\ref{r:3.2}.

\emph{(ii)} Being bipartite, the graph has no odd cycle. A $4$-cycle would be $X\,A\,X'\,A'$ with $X\ne X'$ of size $k-1$ contained in both $A$ and $A'$; then $|X\cup X'|=k$ and $A=A'=X\cup X'$, a contradiction. To construct a $6$-cycle, choose distinct leaves $a,b,c$ and $W\subseteq[n]\setminus\{a,b,c\}$ with $|W|=k-2$, which is possible since $k-2\le n-3$, and take $W\cup\{a\}$, $W\cup\{a,b\}$, $W\cup\{b\}$, $W\cup\{b,c\}$, $W\cup\{c\}$, $W\cup\{a,c\}$. Consecutive sets in this list (cyclically) are comparable and alternate in size, so the six vertices form a $6$-cycle in $F_k(S_n)$.

\emph{(iii)} Every $\sigma\in\Sym(n)$ acts on the subsets of $[n]$ and preserves containment, so it induces an automorphism of $B(n;k-1,k)$; we use these to move one edge onto another. An edge of $B(n;k-1,k)$ is a pair $\{X,A\}$ with $X\subset A$, $|X|=k-1$ and $|A|=k$, so it is determined by $X$ together with the single element $x$ of $A\setminus X$. Let $\{X,A\}$ and $\{X',A'\}$ be two edges, and write $A=X\cup\{x\}$ and $A'=X'\cup\{x'\}$. Since $|X|=|X'|=k-1$, $x\notin X$ and $x'\notin X'$, there is a permutation $\sigma\in\Sym(n)$ with $\sigma(X)=X'$ and $\sigma(x)=x'$, obtained by mapping $X$ onto $X'$ bijectively, sending $x$ to $x'$, and matching the remaining $n-k$ elements of $[n]\setminus(X\cup\{x\})$ with those of $[n]\setminus(X'\cup\{x'\})$. Then $\sigma(A)=\sigma(X)\cup\{\sigma(x)\}=A'$, so the induced automorphism carries the edge $\{X,A\}$ to the edge $\{X',A'\}$. Hence $F_k(S_n)$ is edge-transitive. Moreover, $F_k(S_n)$ is connected, since by Lemma~\ref{r:3.6} any two vertices are at finite distance.  Thus, Watkins's theorem~\cite{b55} gives $\kappa(F)=\lambda(F)=\delta(F)$, which is $k$ by Corollary~\ref{r:3.2}.

\emph{(iv)} Each side of the bipartition is independent, so $\alpha\ge\binom nk$, and $\nu\le\binom{n}{k-1}$, the size of the smaller side. For the reverse, let $\cs{A}$ be a family of $a$ many $(k-1)$-subsets. The number of edges leaving $\cs{A}$ is equal to $(n-k+1)a$ and each neighbor absorbs at most $k$ of them, so $|N(\cs{A})|\ge(n-k+1)a/k\ge a$ because $n-k+1\ge k$. Thus, Hall's condition holds, so there is a matching saturating the smaller side and $\nu(F)=\binom{n}{k-1}$; K\H{o}nig's theorem~\cite{b41} gives $\tau(F)=\nu(F)$, and the Gallai identity $\alpha(F)+\tau(F)=|V(F)|$ together with $|V(F)|=\binom{n}{k-1}+\binom nk$ gives $\alpha(F)=\binom nk$.

\emph{(v)} In the range $2\le k\le(n+1)/2$ the case $n\le4$ leaves only $k=2$ with $n\in\{3,4\}$, and $F_2(S_3)=C_6$ and $F_2(S_4)$ are planar by Proposition~\ref{r:3.5}. For $n\ge5$ we exhibit a nonplanar subgraph. If $k=2$, Corollary~\ref{r:3.4} identifies $F_2(S_n)$ with the total subdivision of $K_n$, which is nonplanar because $K_n$ is. Suppose $k\ge3$ and $n\ge k+3$. Fix a set $L$ of five leaves of $S_n$ and a set $W$ of $k-2$ leaves disjoint from $L$, and consider the configurations that keep $k-2$ tokens fixed at the leaves of $W$ and the remaining two tokens on $\{0\}\cup L$. Two of them differ by an elementary move precisely when the two free tokens do, so they induce a copy of $F_2(S_5)$, the total subdivision of $K_5$, which is nonplanar. It remains to check the pairs $(k,n)$ of the range that the previous paragraph misses, namely those with $k\ge3$ and $n<k+3$, where $W$ and $L$ do not fit. In the range $2k\le n+1$ we have $n\ge2k-1$, and $2k-1\ge k+3$ holds exactly when $k\ge4$; so for $k\ge4$ there is nothing left to check, and the only remaining pair is $k=3$ with $n=5$. There $F_3(S_5)\cong B(5;2,3)$ contains a subdivision of $K_{3,3}$, whose branch vertices are the three pairs contained in $\{1,2,3\}$ and the three triples containing $\{1,2\}$, that is, $\{\{1,2\},\{1,3\},\{2,3\}\}$ is one side and $\{\{1,2,3\},\{1,2,4\},\{1,2,5\}\}$ is the other side. Five of the nine required connections are edges, since $\{1,2\}$ is contained in each of the three triples, and $\{1,3\}$ and $\{2,3\}$ are contained in $\{1,2,3\}$. The remaining four are the paths
\[
\{u,3\}\;-\;\{u,3,x\}\;-\;\{u,x\}\;-\;\{1,2,x\},\qquad u\in\{1,2\},\ x\in\{4,5\},
\]
whose eight interior vertices $\{1,3,4\},\{1,4\},\{1,3,5\},\{1,5\},\{2,3,4\},\{2,4\},\{2,3,5\},\{2,5\}$ are pairwise distinct and are not branch vertices. Hence $F_3(S_5)$ is not planar either.
\end{proof}

\subsection{Automorphisms}\label{sec:A4}

The automorphism group of $F_k(S_n)$ is known. The case $k=2$, $n\ge4$ was obtained by Ibarra and Rivera~\cite{b34}, and in general $\Aut(F_k(G))$ was determined in~\cite{b24} for every connected $(C_4,\text{diamond})$-free graph $G$, where it is $\Aut(G)$ when $2k\ne|V(G)|$ and $\Aut(G)\times\mathbb{Z}_2$ otherwise; the star is a tree, hence $(C_4,\text{diamond})$-free, and $|V(S_n)|=n+1$, so the next proposition is the case $G=S_n$ of that result. The proof below is specific to the star: the identification with $B(n;k-1,k)$ reduces the question to the automorphism group of the Johnson graph $J(n,k-1)$, whereas the argument of~\cite{b24} reconstructs $G$ from $F_k(G)$ for the whole class.

\begin{proposicion}[see~\cite{b24}]\label{r:3.11}
For $2\le k\le(n+1)/2$,
\[
\Aut\bigl(F_k(S_n)\bigr)\;\cong\;
\begin{cases}
\mathbb{Z}_2\times\Sym(n), & \text{if } 2k=n+1,\\
\Sym(n), & \text{otherwise.}
\end{cases}
\]
\end{proposicion}

\begin{proof}
One inclusion is known in general. Every $f\in\Aut(G)$ induces an automorphism $\iota(f)$ of $F_k(G)$ by acting on subsets, and $\iota$ is an injective group homomorphism~\cite{b34}; moreover the complementation $\mathfrak{c}$ of Section~\ref{sec:2} is an automorphism of $F_k(G)$ precisely when $2k=|V(G)|$, and in that case it is not of the form $\iota(f)$ and commutes with every such map~\cite{b24}. Since $\Aut(S_n)\cong\Sym(n)$ for $n\ge3$ and $|V(S_n)|=n+1$, this gives $\Sym(n)\le\Aut(F_k(S_n))$ in general, and $\mathbb{Z}_2\times\Sym(n)\le\Aut(F_k(S_n))$ when $2k=n+1$. It remains to prove the reverse inclusions, and we do so in $B(n;k-1,k)$.

For the reverse inclusion, suppose first $2k\ne n+1$. Then the two sides have different degrees, so every automorphism $\varphi$ preserves each side. Its restriction to the $(k-1)$-side is a permutation preserving the relation ``$X$ and $X'$ lie in a common $k$-set'', which holds precisely when $|X\cap X'|=k-2$; that is, $\varphi$ restricts to an automorphism of the Johnson graph $J(n,k-1)$. Now $\Aut(J(n,m))\cong\Sym(n)$ except when $n=2m$ (see Ramras and Donovan~\cite{bRD}), and $n=2(k-1)$ would give $2k=n+2$, which our range $2k\le n+1$ excludes. So $\varphi$ agrees on the $(k-1)$-side with some $\sigma\in\Sym(n)$, and since $\varphi$ preserves containment it agrees with $\sigma$ on the $k$-side as well. Hence $\Aut\cong\Sym(n)$.

If $2k=n+1$, the two sides have the same size and degree and an automorphism may interchange them; but for any $\varphi$ either $\varphi$ or $\mathfrak{c}\circ\varphi$ preserves the sides, and here $n=2k-1\ne2(k-1)$, so the previous paragraph applies to it. Therefore $\Aut=\langle\mathfrak{c}\rangle\times\Sym(n)\cong\mathbb{Z}_2\times\Sym(n)$.
\end{proof}

For $n=3$ we have $2k=n+1$ and $F_2(S_3)\cong C_6$, so $\Aut(F_2(S_3))\cong\mathbb{Z}_2\times\Sym(3)$ is the dihedral group of $C_6$, as the proposition says.

\subsection{Spectrum}\label{sec:A5}

The \emph{spectrum} of $G$ is the multiset of eigenvalues of its adjacency matrix $A(G)$; the Laplacian matrix $L(G)$ and the algebraic connectivity $a(G)$ were recalled in Subsection~\ref{sec:3.3}.

The adjacency spectrum of $F_k(S_n)$ follows from the singular values of the inclusion matrix of $(k-1)$-sets against $k$-sets, which are classical (see, e.g.,~\cite[Section~9.1]{b11}), and the Laplacian spectrum then follows because $F_k(S_n)$ is bipartite and biregular. We include both derivations: in the inclusion model they are short, since the relevant Gram matrix is a multiple of the identity plus the adjacency matrix of a Johnson graph, and they give the multiplicities for every $n$ and $k$, which is what Proposition~\ref{r:3.13} and Corollary~\ref{r:sptrees} need.

\begin{proposicion}\label{r:3.12}
For $2\le k\le(n+1)/2$, the spectrum of $F_k(S_n)$ consists of
\[
\pm\sqrt{(k-j)(n-k+1-j)}\quad\text{with multiplicity}\quad\binom nj-\binom{n}{j-1},\qquad 0\le j\le k-1,
\]
together with $0$ with multiplicity $\binom nk-\binom{n}{k-1}$.
\end{proposicion}
\begin{proof}
Let $W$ be the inclusion matrix of the $(k-1)$-subsets against the $k$-subsets of $[n]$, so that, ordering the vertices by side, the adjacency matrix of $B(n;k-1,k)$ is $M=\left(\begin{smallmatrix}0&W\\W^{\top}&0\end{smallmatrix}\right)$. Write $N=n-k+1$ and $m=\binom{n}{k-1}$.

We first compute $WW^{\top}$. For $(k-1)$-subsets $X,X'$ the entry $(WW^{\top})_{X,X'}$ counts the $k$-subsets containing $X\cup X'$; there are $N$ of them if $X=X'$, exactly one if $|X\cap X'|=k-2$, and none otherwise. Since two distinct $(k-1)$-sets are adjacent in $J(n,k-1)$ precisely when they meet in $k-2$ points,
\[
WW^{\top}\;=\;N\,I\;+\;A\bigl(J(n,k-1)\bigr).
\]
The eigenvalues of $J(n,k-1)$ are $\theta_j=(k-1-j)(N-j)-j$, with multiplicity $\binom nj-\binom{n}{j-1}$, for $0\le j\le\min\{k-1,N\}$ (see~\cite[Section~9.1]{b11}), and $\min\{k-1,N\}=k-1$ because $k\le(n+1)/2$. Since $WW^{\top}$ is $N\,I$ plus $A(J(n,k-1))$, it has the same eigenvectors as $A(J(n,k-1))$, with eigenvalues shifted by $N$; as $(k-1-j)(N-j)=(k-j)(N-j)-(N-j)$, these are
\[
N+\theta_j\;=\;(k-j)(N-j),\qquad 0\le j\le k-1,
\]
with the same multiplicities, and they are all positive, since $N-j\ge n-2k+2\ge1$.

Let $x_1,\dots,x_m$ be an orthonormal basis of eigenvectors of $WW^{\top}$, with $WW^{\top}x_i=\sigma_i^2x_i$ and $\sigma_i>0$, and put $y_i=W^{\top}x_i/\sigma_i$. Then $Wy_i=\sigma_ix_i$ and $W^{\top}x_i=\sigma_iy_i$, 
and $y_1,\dots,y_m$ are orthonormal, since $\langle y_i,y_l\rangle=\langle x_i,WW^{\top}x_l\rangle/(\sigma_i\sigma_l)=\sigma_l\langle x_i,x_l\rangle/\sigma_i$, which is $1$ if $i=l$ and $0$ otherwise.
Therefore $M(x_i,\pm y_i)=\pm\sigma_i\,(x_i,\pm y_i)$. Moreover, if $v$ is orthogonal to $y_1,\dots,y_m$, then $\langle x_i,Wv\rangle=\sigma_i\langle y_i,v\rangle=0$ for every $i$, so $Wv=0$ and $M(0,v)=0$; these vectors form a space of dimension $\binom nk-m$. All these vectors are pairwise orthogonal, and $2m+\bigl(\binom nk-m\bigr)=|V(F_k(S_n))|$, so we have found all the eigenvalues: $\pm\sqrt{(k-j)(N-j)}$ with multiplicity $\binom nj-\binom{n}{j-1}$, and $0$ with multiplicity $\binom nk-\binom{n}{k-1}$.
\end{proof}

For $k=2$ the proposition gives $\pm\sqrt{2(n-1)}$ once, $\pm\sqrt{n-2}$ with multiplicity $n-1$, and $0$ with multiplicity $\binom n2-n$, which is the spectrum of the total subdivision of $K_n$ (Corollary~\ref{r:3.4}; see~\cite{b16}). The Laplacian spectrum (Proposition~\ref{r:3.13}) follows from the same singular value decomposition.

\begin{proof}[Proof of Proposition~\ref{r:3.13}]
We keep the notation of the proof of Proposition~\ref{r:3.12}. The centered vertices have degree $N=n-k+1$ and the center-free ones have degree $k$, so, ordering the vertices by side, the Laplacian is $L=\left(\begin{smallmatrix}NI&-W\\-W^{\top}&kI\end{smallmatrix}\right)$. For each $i$,
\[
L\begin{pmatrix}x_i\\0\end{pmatrix}=\begin{pmatrix}Nx_i\\-\sigma_i y_i\end{pmatrix},\qquad L\begin{pmatrix}0\\y_i\end{pmatrix}=\begin{pmatrix}-\sigma_i x_i\\ky_i\end{pmatrix},
\]
so the plane $P_i$ spanned by $(x_i,0)$ and $(0,y_i)$ is $L$-invariant, and $L$ acts on it as $\left(\begin{smallmatrix}N&-\sigma_i\\-\sigma_i&k\end{smallmatrix}\right)$. If $\sigma_i^2=(k-j)(N-j)$, this matrix has trace $N+k=n+1$ and determinant $Nk-(k-j)(N-j)=j(n+1-j)$, so its eigenvalues are $j$ and $n+1-j$. On the space of vectors $(0,v)$ with $v$ orthogonal to $y_1,\dots,y_m$, which satisfy $Wv=0$, $L$ acts as multiplication by $k$. The planes $P_1,\dots,P_m$ and this space are pairwise orthogonal, because the $x_i$ are orthonormal and so are the $y_i$, and their dimensions add up to $|V(F_k(S_n))|$. Hence the Laplacian eigenvalues are $j$ and $n+1-j$, each with multiplicity $\binom nj-\binom{n}{j-1}$, for $0\le j\le k-1$, and $k$, with multiplicity $\binom nk-\binom{n}{k-1}$.

The index $j=0$ gives the eigenvalues $0$ and $n+1$, each once, and, since $k\ge2$, the index $j=1$ gives the eigenvalue $1$; all the other eigenvalues are at least $2$. Hence $a(F_k(S_n))=1$, and by~\eqref{eq:1} this holds for every $2\le k\le n-1$.
\end{proof}

The equality $a(F_k(S_n))=a(S_n)=1$ is an instance of a general fact: $a(F_k(G))=a(G)$ for every graph $G$ on $n$ vertices and $1\le k\le n-1$. This is a consequence of Aldous' spectral gap conjecture, proved by Caputo, Liggett and Richthammer~\cite{b13}; see~\cite{b17} for this connection and for an algebraic proof in the case of trees~\cite[Theorem~3.6]{b17}.

\section{Known results from design and coding theory}\label{sec:B}

We collect here, in the language of $F_k(S_n)$, the results from design and coding theory that Theorem~\ref{r:1.1} makes available for the packing number. 

By Theorem~\ref{r:1.1}, $\rho(F_k(S_n))=\mathcal{D}(k-1,k,n+1)=A(n+1,4,k)$, and the Steiner bound of Proposition~\ref{r:5.1} is the first consequence. The bound of Proposition~\ref{r:5.1} can be improved by applying the same count point by point. This is the Johnson bound~\cite{b37}, which we state in the form we will use.

\begin{proposicion}[Johnson bound]\label{r:5.2}
For $v\ge k>t\ge 1$ we have
\[
\mathcal{D}(t,k,v)\le\left\lfloor \frac{v}{k}\,\mathcal{D}(t-1,k-1,v-1)\right\rfloor,
\]
and $\mathcal{D}(0,k,v)=1$. Consequently, for $2\le k\le n$,
\[
\rho(F_k(S_n))\le J(n+1,k):=\left\lfloor\frac{n+1}{k}\left\lfloor\frac{n}{k-1}\left\lfloor\cdots\left\lfloor\frac{n-k+3}{2}\right\rfloor\cdots\right\rfloor\right\rfloor\right\rfloor.
\]
\end{proposicion}

\begin{proof}
Since the empty set is contained in every block, a partial $S(0,k,v)$ has only one block. Let $\mathcal{B}$ be a partial $S(t,k,v)$ with $t\ge 1$ on a set $X$ of $v$ points and, for $x\in X$, let $r_x$ be the number of blocks containing $x$. The family $\{B\setminus\{x\}: x\in B\in\mathcal{B}\}$ is a partial $S(t-1,k-1,v-1)$ on $X\setminus\{x\}$, since two of its members sharing a $(t-1)$-set $T$ would come from two blocks sharing $T\cup\{x\}$. Hence $r_x\le\mathcal{D}(t-1,k-1,v-1)$, and counting the pairs $(x,B)$ with $x\in B$ gives \[k|\mathcal{B}|=\sum_{x\in X}r_x\le v\,\mathcal{D}(t-1,k-1,v-1).\] The formula for $\rho(F_k(S_n))$ is obtained by iterating from $t=k-1$ down to $t=0$ and using Theorem~\ref{r:1.1}.
\end{proof}

The Johnson bound is not always attained; for instance, $J(10,6)=41$, whereas $\rho(F_6(S_9))=A(10,4,6)=A(10,4,4)=30$~\cite{b12}. For $k=3$ and $k=4$ the exact value of $\mathcal{D}(k-1,k,v)$ is known for all $v$, and for $k\ge 5$ the problem is open in general, as we now describe.

For $k=3$ the exact value of $\rho(F_3(S_n))$ is known for every $n$. Kirkman~\cite{b39} proved that a Steiner triple system $S(2,3,v)$ exists if and only if $v\equiv 1,3\pmod{6}$. For the remaining values of $v$, the maximum packing was determined by Schönheim~\cite{b52} and Spencer~\cite{b53}. Combining this with Theorem~\ref{r:1.1} we obtain:

\begin{teorema}[\cite{b39,b52,b53}]\label{r:5.11}
Let $n\ge 3$ and $v=n+1$. Then
\[
\rho(F_3(S_n))=\mathcal{D}(2,3,v)=
\begin{cases}
\left\lfloor\dfrac{v}{3}\left\lfloor\dfrac{v-1}{2}\right\rfloor\right\rfloor-1, & \text{if } v\equiv 5\pmod{6},\\[3ex]
\left\lfloor\dfrac{v}{3}\left\lfloor\dfrac{v-1}{2}\right\rfloor\right\rfloor, & \text{otherwise.}
\end{cases}
\]
In particular, $\rho(F_3(S_n))=\frac{1}{3}\binom{n+1}{2}$ if and only if $n\equiv 0,2\pmod{6}$.
\end{teorema}

Thus, for $k=3$ the Johnson bound (Proposition~\ref{r:5.2}) is attained except when $n+1\equiv 5\pmod{6}$, in which case the value is one less; for instance $\rho(F_3(S_4))=2$, although $J(5,3)=3$.

For $k=4$, Hanani~\cite{b32} proved that a Steiner quadruple system $S(3,4,v)$ exists if and only if $v\equiv 2,4\pmod{6}$. The value of $\mathcal{D}(3,4,v)$ for all $v$ was determined after a series of works, among them that of Bitan and Etzion~\cite{b9}; the last remaining cases were settled by Bao and Ji~\cite{b7}, where the earlier references can be found. The result is that the Johnson bound is attained after a single step of the recursion of Proposition~\ref{r:5.2}.

\begin{teorema}[\cite{b32,b9,b7}]\label{r:5.15}
Let $n\ge 4$ and $v=n+1$. Then $\mathcal{D}(3,4,v)=\bigl\lfloor\frac{v}{4}\,\mathcal{D}(2,3,v-1)\bigr\rfloor$, that is,
\[
\rho(F_4(S_n))=\left\lfloor\frac{n+1}{4}\,\rho(F_3(S_{n-1}))\right\rfloor.
\]
Explicitly,
\[
\rho(F_4(S_n))=
\begin{cases}
\left\lfloor\dfrac{v}{4}\left(\left\lfloor\dfrac{v-1}{3}\left\lfloor\dfrac{v-2}{2}\right\rfloor\right\rfloor-1\right)\right\rfloor, & \text{if } v\equiv 0\pmod{6},\\[3ex]
\left\lfloor\dfrac{v}{4}\left\lfloor\dfrac{v-1}{3}\left\lfloor\dfrac{v-2}{2}\right\rfloor\right\rfloor\right\rfloor, & \text{otherwise.}
\end{cases}
\]
In particular, $\rho(F_4(S_n))=\frac{1}{4}\binom{n+1}{3}$ if and only if $n\equiv 1,3\pmod{6}$.
\end{teorema}

The explicit form follows from the recursive one and from Theorem~\ref{r:5.11}, since for $v\equiv 0\pmod{6}$ we have $v-1\equiv 5\pmod{6}$ and $\mathcal{D}(2,3,v-1)$ is one less than the Johnson bound, while in the remaining cases it coincides with it. For instance, for $n=11$ the Johnson bound is $J(12,4)=54$, whereas $\rho(F_4(S_{11}))=\lfloor 3\cdot 17\rfloor=51$.

For $k\ge 5$ no formula for $\mathcal{D}(k-1,k,v)$ is known. For small values, the best known bounds for $A(v,4,k)$ are in the tables of constant-weight codes~\cite{b12,b10,b1}. Östergård~\cite{b45} classified the optimal constant-weight codes of small length; in particular, $\mathcal{D}(4,5,v)$ is known exactly for $v\le 12$, and $\rho(F_5(S_{12}))=\mathcal{D}(4,5,13)$, with $123\le\mathcal{D}(4,5,13)\le 129$ (the lower bound from~\cite{b10}, the upper bound from~\cite{b45}), is the first unknown value of the family $\rho(F_5(S_n))$. Table~\ref{tab:1} collects the known values and bounds for $k\in\{5,6\}$ and $10\le v\le 16$, together with the Johnson bound of Proposition~\ref{r:5.2}. The values coming from Steiner systems are, by Proposition~\ref{r:5.1}, $\rho(F_5(S_{10}))=66$, given by the Witt design $S(4,5,11)$, and $\rho(F_6(S_{11}))=132$, given by $S(5,6,12)$ (Example~\ref{r:5.3}(d)); in both the Johnson bound is attained, and in the unknown range it is far from the best known lower bounds.

\begin{table}[t]
\centering
\begin{tabular}{ll|ccccccc}
\toprule
 & $v=n+1$ & $10$ & $11$ & $12$ & $13$ & $14$ & $15$ & $16$\\
\midrule
$k=5$ & $\rho(F_5(S_n))=A(v,4,5)$ & $36$ & $66$ & $80$ & $123$--$129$ & $\ge 169$ & $\ge 243$ & $\ge 322$\\
 & Johnson bound $J(v,5)$ & $\le 36$ & $\le 66$ & $\le 84$ & $\le 140$ & $\le 182$ & $\le 273$ & $\le 336$\\
\midrule
$k=6$ & $\rho(F_6(S_n))=A(v,4,6)$ & $30$ & $66$ & $132$ & $\ge 166$ & $\ge 278$ & $\ge 399$ & $\ge 624$\\
 & Johnson bound $J(v,6)$ & $\le 41$ & $\le 66$ & $\le 132$ & $\le 182$ & $\le 326$ & $\le 455$ & $\le 728$\\
\bottomrule
\end{tabular}
\caption{Values of $\rho(F_k(S_n))=A(n+1,4,k)$ for $k=5,6$ and $10\le n+1\le 16$. Values without ``$\ge$'' or an interval are exact; those marked with ``$\ge$'' are the best known lower bounds~\cite{b10}, and the interval gives the best known lower~\cite{b10} and upper~\cite{b45} bounds. The second row of each block is the Johnson upper bound of Proposition~\ref{r:5.2}.}
\label{tab:1}
\end{table}

\section{Omitted proofs and examples of packings}\label{sec:C}

This appendix contains the proofs of Proposition~\ref{r:6.5} and Lemma~\ref{r:6.2}, examples of the constructions of Subsection~\ref{sec:5.2}, and, for $k=3$, their comparison with the exact value of the packing number. We begin with the divisibility conditions of Proposition~\ref{r:6.5}.

\begin{proof}[Proof of Proposition~\ref{r:6.5}]
(i) Let $u=v-k+1\ge1$ and $a=k-i$. As $i$ ranges over $0,\dots,k-2$, $a$ ranges over $2,\dots,k$, and $v-i=u+a-1$ and $k-1-i=a-1$. The $i$-th condition is then $a\mid\binom{u+a-1}{a-1}$, and since
\[
\binom{u+a-1}{a-1}=\frac{(u+1)(u+2)\cdots(u+a-1)}{(a-1)!},
\]
it is equivalent to
\begin{equation}\label{eq:13}
a!\;\Big|\;\prod_{j=1}^{a-1}(u+j).
\end{equation}
Suppose first that $\gcd(u,P_k)=1$, and let $2\le a\le k$. The product of the $a$ consecutive integers $u,u+1,\dots,u+a-1$ is $a!\binom{u+a-1}{a}$, so $a!$ divides $u\prod_{j=1}^{a-1}(u+j)$. Every prime dividing $a!$ is at most $k$, so it divides $P_k$ and not $u$; hence $\gcd(u,a!)=1$, and~\eqref{eq:13} follows. Conversely, suppose that a prime $p\le k$ divides $u$, and take $a=p$. The integers $u+1,\dots,u+p-1$ are congruent to $1,\dots,p-1$ modulo $p$, so $p$ does not divide their product, and~\eqref{eq:13} fails, since $p\mid p!$.

(ii) By \eqref{eq:10}, the divisibility conditions of an $S(k-2,k-1,2k-1)$ are
\[
(k-1-i)\;\Big|\;\binom{2k-1-i}{k-2-i},\qquad 0\le i\le k-3.
\]
With $j=i+1$ they read $(k-j)\mid\binom{2k-j}{k-1-j}$ for $1\le j\le k-2$, which are the conditions~\eqref{eq:adm} for $v=2k$ with $j\ge1$. The remaining condition, $j=0$, says $k\mid\binom{2k}{k-1}$, and it always holds, since $\binom{2k}{k-1}=\frac{k}{k+1}\binom{2k}{k}=k\,C_k$. Thus, the divisibility conditions of $S(k-2,k-1,2k-1)$ are equivalent to $2k$ being $k$-admissible, and by (i) this is equivalent to $\gcd(k+1,P_k)=1$, that is, to no prime $p\le k$ dividing $k+1$. Since $k+1\ge3$, this happens if and only if $k+1$ is prime, because a composite $k+1$ has a prime factor $p\le\sqrt{k+1}\le k$.
\end{proof}

Next we prove Lemma~\ref{r:6.2}. As in Subsection~\ref{sec:6.2}, $[n]=[2k-1]$ denotes the set of leaves and $\lambda_i$ the numbers in~\eqref{eq:10}.

\begin{proof}[Proof of Lemma~\ref{r:6.2}]
(i) Take $i=k-3$ in \eqref{eq:10}, which is valid because $k\ge3$. With $t=k-2$, $r=k-1$ and $v=2k-1$ we obtain
\[
\lambda_{k-3}=\binom{(2k-1)-(k-3)}{1}\Big/\binom{(k-1)-(k-3)}{1}=\frac{k+2}{2},
\]
which must be an integer; hence $k$ is even.

(ii) Two distinct $(k-1)$-subsets of a $k$-set $A$ intersect in $k-2$ elements, so they share a $(k-2)$-subset. If both were blocks, that $(k-2)$-set would lie in two blocks, which contradicts the Steiner property.

(iii) For a $(k-1)$-subset $X$ we write $g(X)=\bigl|\{Y\in\mathcal{S}:Y\cap X=\varnothing\}\bigr|$. Every block disjoint from $X$ is contained in $[n]\setminus X$, which is a $k$-set, so by (ii) there is at most one, implying $g(X)\in\{0,1\}$. By inclusion--exclusion over the points of $X$,
\[
g(X)=\sum_{J\subseteq X}(-1)^{|J|}\,\bigl|\{Y\in\mathcal{S}:J\subseteq Y\}\bigr|.
\]
If $|J|=j\le k-2$, this number is $\lambda_j$, which depends only on $j$; and if $J=X$, since $|X|=k-1=|Y|$, the number is $1$ if $X\in\mathcal{S}$ and $0$ otherwise. Therefore
\[
g(X)=c+(-1)^{k-1}\,[\,X\in\mathcal{S}\,],\qquad c:=\sum_{j=0}^{k-2}(-1)^j\binom{k-1}{j}\lambda_j,
\]
where $c$ depends only on $k$ (through the numbers $\lambda_j$), and not on $X$. By (i), $k$ is even, so $(-1)^{k-1}=-1$ and $g(X)=c-[\,X\in\mathcal{S}\,]$.

We now bound $c$ from both sides. There exists some $(k-1)$-set that is not a block, since
\[
\binom{2k-1}{k-1}=\frac{k+1}{k-1}\binom{2k-1}{k-2}=(k+1)\,|\mathcal{S}|>|\mathcal{S}|;
\]
for such an $X$ we have $c=g(X)\in\{0,1\}$, so $c\le1$. And taking $X\in\mathcal{S}$, which is nonempty, we obtain $c-1=g(X)\ge0$, so $c\ge1$. Therefore $c=1$, and consequently $g(X)=1$ if $X\notin\mathcal{S}$, and $g(X)=0$ if $X\in\mathcal{S}$, which is (iii).

Finally, we prove (iv). If $A$ is a $k$-set and $X=[n]\setminus A$, which is a $(k-1)$-set, then $A$ contains a block if and only if there exists a block disjoint from $X$, that is, if and only if $g(X)=1$, which by (iii) happens exactly when $X\notin\mathcal{S}$.
\end{proof}

The following example compares the two constructions of Subsection~\ref{sec:5.2}. In~(a) the size of a sum class depends on the group, and in (b)--(d) one of the two constructions is optimal and the other is not.

\begin{ejemplo}\label{r:5.10}\leavevmode
\begin{enumerate}[label=(\alph*)]
\item $k=3$, $v=9$. Here $|\mathcal{P}_0(\mathbb{Z}_9,3)|=10$, whereas $|\mathcal{P}_0(\mathbb{Z}_3^2,3)|=12=\frac{1}{3}\binom{9}{2}$, since three distinct points of $\mathbb{Z}_3^2$ sum to $0$ if and only if they are collinear, so that $\mathcal{P}_0(\mathbb{Z}_3^2,3)$ is the affine plane of order $3$, an $S(2,3,9)$. Thus $\rho(F_3(S_8))=12$ is attained by a sum class, but only if the group is chosen well.
\item $k=3$, $v=8$. Here $K(8,3)=7<8=\rho(F_3(S_7))$ (Example~\ref{r:5.3}(c)), so no sum class is a maximum packing of $F_3(S_7)$. The class $\mathcal{P}_0(\mathbb{Z}_7,3)$ has $5$ blocks, and its leave consists of the pairs $\{x,-2x\}$, $x\neq 0$, which form the $6$-cycle $1,5,4,6,2,3$; the three blocks $\{1,5,\infty\},\{4,6,\infty\},\{2,3,\infty\}$ through a new point $\infty$ complete it to a $\mathbb{Z}_7+1$ point construction with $8$ blocks, which is optimal.
\item $k=4$, $v=8$. Since $\gcd(8,4)=4$, $\mathcal{R}=\{1,2\}$ and $w_4(2)=2^3-1=7$, so $K(8,4)=\frac{1}{8}\bigl(\binom{8}{4}+\binom{4}{2}\cdot 7\bigr)=14=\frac{1}{4}\binom{8}{3}$, attained by $\mathcal{P}_0(\mathbb{Z}_2^3,4)$. Indeed, four distinct points of $\mathbb{Z}_2^3$ sum to $0$ if and only if they form an affine plane, so this class consists of the $14$ planes of $AG(3,2)$, an $S(3,4,8)$, and $\rho(F_4(S_7))=14$. With $G=\mathbb{Z}_7$ and one new point one obtains only $11$ blocks, and no choice of $c$ does better.
\item $k=4$, $v=9$. Since $\gcd(9,4)=1$, every sum class of a group of order $9$ has $\frac{1}{9}\binom{9}{4}=14$ blocks, whereas $\rho(F_4(S_8))=18$ (Theorem~\ref{r:5.15}). The class $\mathcal{P}_0(\mathbb{Z}_8,4)$ has $9$ blocks,
\begin{gather*}
    \{0,1,2,5\},\{0,1,3,4\},\{0,3,6,7\},\{0,4,5,7\},\{1,2,6,7\},\\ \{1,3,5,7\},\{1,4,5,6\},\{2,3,4,7\},\{2,3,5,6\},
\end{gather*}
which cover $36$ of the $56$ triples of $\mathbb{Z}_8$; the block $\{0,2,4,6\}$, whose four triples lie in the leave, and the eight blocks $X\cup\{\infty\}$ with $X\in\{016,027,035,123,147,245,346,567\}$ give a $\mathbb{Z}_8+1$ point construction with $18$ blocks, which is optimal.
\end{enumerate}
\end{ejemplo}

Finally, for $k=3$ the leave and the size of a sum class can be described explicitly, which allows us to compare the sum classes with the exact value.

\begin{lema}\label{r:5.12}
Let $\Gamma$ be an abelian group of order $v$, $c\in \Gamma$ and $T_c=|\{x\in \Gamma:\ 3x=c\}|$.
\begin{enumerate}[label=(\roman*)]
\item A pair $\{x,y\}\subseteq \Gamma$ is not contained in any block of $\mathcal{P}_c(\Gamma,3)$ if and only if $y=c-2x$ or $x=c-2y$.
\item $|\mathcal{P}_c(\Gamma,3)|=\dfrac{v^2-3v+2T_c}{6}$.
\end{enumerate}
\end{lema}
\begin{proof}
(i) The only possible block containing $\{x,y\}$ is $\{x,y,z\}$ with $z=c-x-y$, and it is a block if and only if $z\notin\{x,y\}$. Moreover, $z=x$ is equivalent to $y=c-2x$, and $z=y$ is equivalent to $x=c-2y$.

(ii) By (i), the uncovered pairs are of the form $\{x,c-2x\}$ with $3x\neq c$. Each of them comes from a single $x$, because if $\{x,c-2x\}=\{y,c-2y\}$ with $y\neq x$, then $y=c-2x$ and $x=c-2y$, thus implying $3x=c$. Therefore, there are exactly $v-T_c$ uncovered pairs, and since each block covers three pairs, $|\mathcal{P}_c(\Gamma,3)|=\bigl(\binom{v}{2}-(v-T_c)\bigr)/3$.
\end{proof}

We include the case $k=3$ of Kl{\o}ve's theorem (Theorem~\ref{r:5.8}) because Lemma~\ref{r:5.12}(ii) gives a direct proof of it and because Proposition~\ref{r:5.14} uses it.

\begin{corolario}\label{r:5.13}
For every $v\ge 4$, with $3^E$ the largest power of $3$ dividing $v$,
\begin{equation}\label{eq:8}
K(v,3)=\frac{v^2-3v+2\cdot 3^E}{6}.
\end{equation}
\end{corolario}
\begin{proof}
By Lemma~\ref{r:5.12}(ii), maximizing $|\mathcal{P}_c(\Gamma,3)|$ amounts to maximizing $T_c$. If $T_c\neq0$, the set $\{x:3x=c\}$ is a coset of the subgroup $H=\{x:3x=0\}$, so $T_c=|H|$. The elements of $H$ have order $1$ or $3$, so $H$ is contained in the Sylow $3$-subgroup of $\Gamma$, which has order $3^E$; hence $T_c\le3^E$. Equality holds for $\Gamma=\Gamma_v^*$ (defined before Theorem~\ref{r:5.8}) and $c=0$, because the Sylow $3$-subgroup of $\Gamma_v^*$ is $\mathbb{Z}_3^E$, all of whose elements satisfy $3x=0$.
\end{proof}

Comparing Corollary~\ref{r:5.13} with the exact values of Theorem~\ref{r:5.11}, we can measure how far the best sum class is from a maximum packing.

\begin{proposicion}\label{r:5.14}
Let $v\ge 4$, let $3^E$ be the largest power of $3$ dividing $v$, and write $\mathcal{D}(v)=\mathcal{D}(2,3,v)$ and $K(v)=K(v,3)$. Then $\mathcal{D}(v)-K(v)$ equals
\[
\begin{array}{ccccccc}
\toprule
v \bmod 6 & 0 & 1 & 2 & 3 & 4 & 5\\
\midrule
\mathcal{D}(v)-K(v) & \frac{v-2\cdot 3^E}{6} & \frac{v-1}{3} & \frac{v-2}{6} & \frac{v-3^E}{3} & \frac{v-4}{6} & \frac{v-5}{3}\\[1ex]
\bottomrule
\end{array}
\]
In particular, $K(v)=\mathcal{D}(v)$ if and only if $v\in\{4,5\}$ or $v=3^a$ or $v=2\cdot 3^a$ with $a\ge 1$, and in every other case $\mathcal{D}(v)-K(v)\ge\frac{v-4}{12}$.
\end{proposicion}

\begin{proof}
By~\eqref{eq:8}, $K(v)=\frac{v^2-3v+2\cdot 3^E}{6}$, with $E=0$ if $3\nmid v$. By Theorem~\ref{r:5.11}, $\mathcal{D}(v)$ equals $\frac{v(v-1)}{6}$ if $v\equiv 1,3$, $\frac{v(v-2)}{6}$ if $v\equiv 0,2$, $\frac{v(v-2)-2}{6}$ if $v\equiv 4$, and $\frac{v(v-1)-8}{6}$ if $v\equiv 5\pmod{6}$. Subtracting, we obtain the table.

Since $v\ge4$, the difference vanishes for $v\equiv0$ if and only if $v=2\cdot3^E$, for $v\equiv3$ if and only if $v=3^E$, and for the other residues if and only if $v\in\{4,5\}$. For the bound, if $v\equiv 0$ and $v\neq 2\cdot 3^E$, then $v/3^E$ is even and greater than $2$, so $v\ge 4\cdot 3^E$ and $\mathcal{D}(v)-K(v)\ge v/12$; if $v\equiv 3$ and $v\neq 3^E$, then $v/3^E\ge 5$ and $\mathcal{D}(v)-K(v)\ge 4v/15$. In the other residues the bound can be read from the table, using $v\ge11$ when $v\equiv5$.
\end{proof}

Thus the difference between the optimum and the best sum class grows linearly with $v$, except for $v\in\{4,5\}$ and for $v$ of the form $3^a$ or $2\cdot 3^a$; for $v=3^a$ the class $\mathcal{P}_0(\mathbb{Z}_3^a,3)$ is the set of lines of the affine geometry $AG(a,3)$, a Steiner triple system, as in Example~\ref{r:5.10}(a). The $\Gamma+j$ points construction does better. Let $\Gamma=\mathbb{Z}_m$ with $m$ odd, so that $T_0\in\{1,3\}$. By Lemma~\ref{r:5.12}(i), the leave of $\mathcal{P}_0(\Gamma,3)$ consists of the pairs $\{x,-2x\}$ with $3x\neq0$, which form the cycles of length at least $3$ of the permutation $x\mapsto-2x$ of $\Gamma\setminus\{0\}$; let $o$ be the number of these cycles of odd length. With one new point $\infty$ we add the blocks $X\cup\{\infty\}$ for $X$ in a maximum matching of each cycle. With two new points $\infty_1,\infty_2$ we add two disjoint such matchings, which together cover every edge of an even cycle and all but one edge of an odd cycle, and the block $\{0,\infty_1,\infty_2\}$. By Lemma~\ref{r:5.12}(ii), the two constructions have $\frac{m^2-T_0}{6}-\frac{o}{2}$ and $\frac{m^2+3m-4T_0+6}{6}-o$ blocks, whereas Theorem~\ref{r:5.11} gives $\mathcal{D}(m+1)=\frac{m^2-T_0}{6}$ and $\mathcal{D}(m+2)=\frac{m^2+3m-4T_0+6}{6}$. So the gap is exactly $o/2$ for $v=m+1$ and exactly $o$ for $v=m+2$.

For $k=4$, the sum class $\mathcal{P}_0(\mathbb{Z}_2^3,4)$ is the Steiner system $AG(3,2)$ of Example~\ref{r:5.10}(c). For $k\ge4$ we have not found a general rule for distributing the leave of a sum class among the new points. A computer search shows that $\mathbb{Z}_7$ with $k-1$ new points yields the Steiner systems $S(k-1,k,k+6)$ for $k=3,4,5,6$, that is, $S(2,3,9)$, $S(3,4,10)$, $S(4,5,11)$ and $S(5,6,12)$; for $k=5$ and $k=6$, however, the sum class contributes only $3$ of the $66$ blocks and $1$ of the $132$ blocks, respectively.

\end{document}